\def\date{}
\catcode`\@=11 
 
\def\nolabels{\def\wrlabel##1{}\def\eqlabel##1{}\def\reflabel##1{}}
\def\writelabels{\def\wrlabel##1{\leavevmode\vadjust{\rlap{\smash%
{\line{{\escapechar=` \hfill\rlap{\sevenrm\hskip.03in\string##1}}}}}}}%
\def\eqlabel##1{{\escapechar-1\rlap{\sevenrm\hskip.05in\string##1}}}%
\def\thlabel##1{{\escapechar-1\rlap{\sevenrm\hskip.05in\string##1}}}%
\def\reflabel##1{\noexpand\llap{\noexpand\sevenrm\string\string\string##1}}}
\nolabels
\global\newcount\secno \global\secno=0
\global\newcount\meqno \global\meqno=1
\global\newcount\mthno \global\mthno=1
\global\newcount\mexno \global\mexno=1
\global\newcount\mquno \global\mquno=1
\global\newcount\tblno \global\tblno=1
\def\newsec#1{\global\advance\secno by1 
\global\subsecno=0\xdef\secsym{\the\secno.}\global\meqno=1\global\mthno=1
\global\mexno=1\global\mquno=1\global\figno=1\global\tblno=1

\bigbreak\medskip\noindent{\bf\the\secno. #1}\writetoca{{\secsym} {#1}}
\par\nobreak\medskip\nobreak}
\xdef\secsym{}
\global\newcount\subsecno \global\subsecno=0
\def\subsec#1{\global\advance\subsecno by1 \global\subsubsecno=0
\xdef\subsecsym{\the\subsecno.}
\bigbreak\noindent{\bf\secsym\the\subsecno. #1}\writetoca{\string\quad
{\secsym\the\subsecno.} {#1}}\par\nobreak\medskip\nobreak}
\xdef\subsecsym{}
\global\newcount\subsubsecno \global\subsubsecno=0
\def\subsubsec#1{\global\advance\subsubsecno by1
\bigbreak\noindent{\it\secsym\the\subsecno.\the\subsubsecno.
                                   #1}\writetoca{\string\quad
{\the\secno.\the\subsecno.\the\subsubsecno.} {#1}}\par\nobreak\medskip\nobreak}
\global\newcount\appsubsecno \global\appsubsecno=0
\def\appsubsec#1{\global\advance\appsubsecno by1 \global\subsubsecno=0
\xdef\appsubsecsym{\the\appsubsecno.}
\bigbreak\noindent{\it\secsym\the\appsubsecno. #1}\writetoca{\string\quad
{\secsym\the\appsubsecno.} {#1}}\par\nobreak\medskip\nobreak}
\xdef\appsubsecsym{}
\def\appendix#1#2{\global\meqno=1\global\mthno=1\global\mexno=1
\global\figno=1\global\tblno=1
\global\subsecno=0\global\subsubsecno=0
\global\appsubsecno=0
\xdef\appname{#1}
\xdef\secsym{\hbox{#1.}}
\bigbreak\bigskip\noindent{\bf Appendix #1. #2}
\writetoca{Appendix {#1.} {#2}}\par\nobreak\medskip\nobreak}
%
%
\def\eqnn#1{\xdef #1{(\secsym\the\meqno)}\writedef{#1\leftbracket#1}%
\global\advance\meqno by1\wrlabel#1}
\def\eqna#1{\xdef #1##1{\hbox{$(\secsym\the\meqno##1)$}}
\writedef{#1\numbersign1\leftbracket#1{\numbersign1}}%
\global\advance\meqno by1\wrlabel{#1$\{\}$}}
\def\eqn#1#2{\xdef #1{(\secsym\the\meqno)}\writedef{#1\leftbracket#1}%
\global\advance\meqno by1$$#2\eqno#1\eqlabel#1$$}
%
%
\def\thm#1{\xdef #1{\secsym\the\mthno}\writedef{#1\leftbracket#1}%
\global\advance\mthno by1\wrlabel#1}
\def\exm#1{\xdef #1{\secsym\the\mexno}\writedef{#1\leftbracket#1}%
\global\advance\mexno by1\wrlabel#1}
%
%
\def\tbl#1{\xdef #1{\secsym\the\tblno}\writedef{#1\leftbracket#1}%
\global\advance\tblno by1\wrlabel#1}
%
\newskip\footskip\footskip14pt plus 1pt minus 1pt 
\def\f@@t{\baselineskip\footskip\bgroup\aftergroup\@foot\let\next}
\setbox\strutbox=\hbox{\vrule height9.5pt depth4.5pt width0pt}
\global\newcount\ftno \global\ftno=0
\def\foot{\global\advance\ftno by1\footnote{$^{\the\ftno}$}}
%
\newwrite\ftfile
\def\footend{\def\foot{\global\advance\ftno by1\chardef\wfile=\ftfile
$^{\the\ftno}$\ifnum\ftno=1\immediate\openout\ftfile=foots.tmp\fi%
\immediate\write\ftfile{\noexpand\smallskip%
\noexpand\item{f\the\ftno:\ }\pctsign}\findarg}%
\def\footatend{\vfill\eject\immediate\closeout\ftfile{\parindent=20pt
\centerline{\bf Footnotes}\nobreak\bigskip\input foots.tmp }}}
\def\footatend{}
%
%
\global\newcount\refno \global\refno=1
\newwrite\rfile
\def\ref{\the\refno\nref}
\def\bref{\nref}
\def\nref#1{\xdef#1{\the\refno}\writedef{#1\leftbracket#1}%
\ifnum\refno=1\immediate\openout\rfile=refs.tmp\fi
\global\advance\refno by1\chardef\wfile=\rfile\immediate
\write\rfile{\noexpand\item{[#1]\ }\reflabel{#1\hskip.31in}\pctsign}\findarg}
\def\findarg#1#{\begingroup\obeylines\newlinechar=`\^^M\pass@rg}
{\obeylines\gdef\pass@rg#1{\writ@line\relax #1^^M\hbox{}^^M}%
\gdef\writ@line#1^^M{\expandafter\toks0\expandafter{\striprel@x #1}%
\edef\next{\the\toks0}\ifx\next\em@rk\let\next=\endgroup\else\ifx\next\empty%
\else\immediate\write\wfile{\the\toks0}\fi\let\next=\writ@line\fi\next\relax}}
\def\striprel@x#1{} \def\em@rk{\hbox{}}
\def\lref{\begingroup\obeylines\lr@f}
\def\lr@f#1#2{\gdef#1{\ref#1{#2}}\endgroup\unskip}

\def\addref#1{\immediate\write\rfile{\noexpand\item{}#1}} 
\def\startrefs#1{\immediate\openout\rfile=refs.tmp\refno=#1}
\def\xref{\expandafter\xr@f}\def\xr@f[#1]{#1}
\def\refs#1{[\r@fs #1{\hbox{}}]}
\def\r@fs#1{\edef\next{#1}\ifx\next\em@rk\def\next{}\else
\ifx\next#1\xref #1\else#1\fi\let\next=\r@fs\fi\next}
%

%
 \newwrite\ffile\global\newcount\figno \global\figno=1
%
%
\def\fig{\the\figno\nfig}
\def\nfig#1{\xdef#1{\secsym\the\figno}%
\writedef{#1\leftbracket \noexpand~\the\figno}%
\ifnum\figno=1\immediate\openout\ffile=figs.tmp\fi\chardef\wfile=\ffile%
\immediate\write\ffile{\noexpand\medskip\noexpand\item{Figure\ \the\figno. }
\reflabel{#1\hskip.55in}\pctsign}\global\advance\figno by1\findarg}
\def\xfig{\expandafter\xf@g}\def\xf@g \penalty\@M\ {}
\def\figs#1{figs.~\f@gs #1{\hbox{}}}
\def\f@gs#1{\edef\next{#1}\ifx\next\em@rk\def\next{}\else
\ifx\next#1\xfig #1\else#1\fi\let\next=\f@gs\fi\next}
%
%
\newwrite\lfile

{\escapechar-1\xdef\pctsign{\string\%}\xdef\leftbracket{\string\{}
\xdef\rightbracket{\string\}}\xdef\numbersign{\string\#}}

\def\writestop{\def\writestoppt{\immediate\write\lfile{\string\pageno%
\the\pageno\string\startrefs\leftbracket\the\refno\rightbracket%
\string\def\string\secsym\leftbracket\secsym\rightbracket%
\string\secno\the\secno\string\meqno\the\meqno}\immediate\closeout\lfile}}
\def\writestoppt{}\def\writedef#1{}

\def\seclab#1{\xdef #1{\the\secno}\writedef{#1\leftbracket#1}\wrlabel{#1=#1}}

\def\subseclab#1{\xdef #1{\secsym\the\subsecno}%
\writedef{#1\leftbracket#1}\wrlabel{#1=#1}}
\def\appsubseclab#1{\xdef #1{\secsym\the\appsubsecno}%
\writedef{#1\leftbracket#1}\wrlabel{#1=#1}}
\def\subsubseclab#1{\xdef #1{\secsym\the\subsecno.\the\subsubsecno}%
\writedef{#1\leftbracket#1}\wrlabel{#1=#1}}
\newwrite\tfile \def\writetoca#1{}
\def\leaderfill{\leaders\hbox to 1em{\hss.\hss}\hfill}
\def\writetoc{\immediate\openout\tfile=toc.tmp
   \def\writetoca##1{{\edef\next{\write\tfile{\noindent ##1
   \string\leaderfill {\noexpand\number\pageno} \par}}\next}}}

%
\catcode`\@=12 
%
%
%
%
%
\def\dbend{{\manual\char127}}
\def\d@nger{\medbreak\begingroup\clubpenalty=10000
    \def\par{\endgraf\endgroup\medbreak} \noindent\hang\hangafter=-2
    \hbox to0pt{\hskip-\hangindent\dbend\hfill}\ninepoint}
\outer\def\danger{\d@nger}

\def\vev#1{\langle #1 \rangle}

\def\darr#1{\raise1.5ex\hbox{$\leftrightarrow$}\mkern-16.5mu #1}

%
%
\def\al{\alpha}
\def\be{\beta}
  \def\Ga{\Gamma}
\def\de{\delta}  \def\De{\Delta}
\def\ep{\epsilon}  
\def\ze{\zeta}

\def\la{\lambda} \def\La{\Lambda}

\def\si{\sigma}

\def\om{\omega}  
%
%

%

%
%

\def\cH{{\cal H}}

\def\cO{{\cal O}}

\def\cW{{\cal W}}

\def\vev#1{\langle #1 \rangle}

\def\proof{\noindent {\it Proof:}\ }
\def\Box{\hbox{$\rlap{$\sqcup$}\sqcap$}}

%

%
%
\def\amsyes{y }

\def\answ{y }

\ifx\answ\amsyes
\input amssym.def


\def\CC{{\Bbb C}}
\def\ZZ{{\Bbb Z}}
\def\NN{{\Bbb N}}
\def\QQ{{\Bbb Q}}

\def\sln{\frak{sl}_N}   \def\hsln{\widehat{\frak{sl}_N}}

\else
\def\ZZ{{Z\!\!\!Z}}              
\def\CC{{I\!\!\!\!C}}
\def\NN{{I\!\!N}}             
\def\QQ{{I\!\!\!\!Q}}

\def\sln{s\ell_n}   \def\hsln{\widehat{s\ell_n}}

\fi
%

%
%

%
%

\newsymbol\ltimes 226E
\newsymbol\rtimes 226F
%
%
%
\def\AdM#1{Adv.\ Math.\ {\bf #1}}

\def\CMP#1{Comm.\ Math.\ Phys.\ {\bf #1}}

\def\IJMP#1{Int.\ J.\ Mod.\ Phys.\ {\bf #1}}

\def\LMP#1{Lett.\ Math.\ Phys.\ {\bf #1}}

\def\NPB#1{Nucl.\ Phys.\ {\bf B#1}}

%

%
%
\def\SMu{\hbox{\lower 3pt\hbox{ \epsffile{su10.eps}}}}
\def\SMs{\hbox{\lower 3pt\hbox{ \epsffile{ss10.eps}}}}
\def\SMd{\hbox{\lower 3pt\hbox{ \epsffile{sd10.eps}}}}

\def\SMS{\leavevmode\vadjust{\rlap{\smash%
{\line{{\escapechar=` \hfill\rlap{\hskip.3in%
                 \hbox{\lower 2pt\hbox{\epsffile{sd10.eps}}}}}}}}}}
\def\SMH{\leavevmode\vadjust{\rlap{\smash%
{\line{{\escapechar=` \hfill\rlap{\hskip.3in%
                 \hbox{\lower 2pt\hbox{\epsffile{su10.eps}}}}}}}}}}
%
%
\def\LW#1{\lower .5pt \hbox{$\scriptstyle #1$}}
\def\LWr#1{\lower 1.5pt \hbox{$\scriptstyle #1$}}
\def\LWrr#1{\lower 2pt \hbox{$\scriptstyle #1$}}
\def\RSr#1{\raise 1pt \hbox{$\scriptstyle #1$}}

%



\def\ABF{1}
\def\AKOSa{2}
\def\AKOSc{3}
\def\AKOSb{4}
\def\BMP{5}
\def\BS{6}
\def\DFJMN{7}
\def\FeFr{8}
\def\FF{9}
\def\FJMOP{10}
\def\FRa{11}
\def\FRb{12}
\def\FRc{13}
\def\FRS{14}
\def\Gar{15}
\def\Gen{16}
\def\HYa{17}
\def\HYb{18}
\def\Jia{19}
\def\Jib{20}
\def\JLMP{21}
\def\LLTa{22}
\def\LLTb{23}
\def\Lit{24}
\def\LPa{25}
\def\LPb{26}
\def\MD{27}
\def\Mila{28}
\def\Milb{29}
\def\Mob{30}
\def\Pol{31}
\def\Sch{32}
\def\SS{33}
\def\SKAO{34}
\def\eqBAa{(2.1)}
\def\eqBAb{(2.2)}
\def\eqBAc{(2.3)}
\def\eqLa{(2.4)}
\def\eqBAd{(2.5)}
\def\eqBAe{(2.6)}
\def\eqBAf{(2.7)}
\def\eqBAfa{(2.8)}
\def\eqBAg{(2.9)}
\def\eqBAh{(2.10)}
\def\thBAa{2.1}
\def\eqBAi{(2.11)}
\def\eqBAj{(2.12)}
\def\eqCaa{(2.13)}
\def\eqCe{(2.14)}
\def\eqCAf{(2.15)}
\def\eqCAfa{(2.16)}
\def\eqCAg{(2.17)}
\def\eqCAh{(2.18)}
\def\thCc{2.2}
\def\KPbasiden{(2.19)}
\def\KPfinid{(2.20)}
\def\KPfandc{(2.21)}
\def\eqKa{(2.22)}
\def\eqBAxa{(2.23)}
\def\KPsing{(2.24)}
\def\eqBBa{(2.25)}
\def\eqBBb{(2.26)}
\def\eqBBaa{(2.27)}
\def\eqBBab{(2.28)}
\def\eqBBac{(2.29)}
\def\eqBBad{(2.30)}
\def\eqBBae{(2.31)}
\def\thBBa{2.3}
\def\eqBBc{(2.32)}
\def\eqBBd{(2.33)}
\def\eqBBdy{(2.34)}
\def\eqBBdz{(2.35)}
\def\eqBBe{(2.36)}
\def\eqBBf{(2.37)}
\def\eqBBg{(2.38)}
\def\eqBBh{(2.39)}
\def\eqBBi{(2.40)}
\def\eqBBj{(2.41)}
\def\eqBBk{(2.42)}
\def\eqBBl{(2.43)}
\def\thBBb{2.4}
\def\eqBBba{(2.44)}
\def\eqBBexch{(2.45)}
\def\eqBBstt{(2.46)}
\def\eqCBaa{(3.1)}
\def\eqCBab{(3.2)}
\def\thCBaa{3.1}
\def\eqCBac{(3.3)}
\def\thCBab{3.2}
\def\eqCBad{(3.4)}
\def\thCBa{3.3}
\def\eqCBa{(3.5)}
\def\eqCBb{(3.6)}
\def\eqCBc{(3.7)}
\def\eqCBd{(3.8)}
\def\eqCBe{(3.9)}
\def\eqCBf{(3.10)}
\def\eqCBg{(3.11)}
\def\eqCBga{(3.12)}
\def\eqCBgb{(3.13)}
\def\eqKb{(3.14)}
\def\eqCBh{(3.15)}
\def\eqKc{(3.16)}
\def\thCBd{3.4}
\def\eqCBh{(3.17)}
\def\eqCBi{(3.18)}
\def\eqCBj{(3.19)}
\def\eqDAa{(4.1)}
\def\eqDAb{(4.2)}
\def\eqDAc{(4.3)}
\def\thDAa{4.1}
\def\eqDAd{(4.4)}
\def\eqDAe{(4.5)}
\def\eqDAf{(4.6)}
\def\thDAb{4.2}
\def\eqDAg{(4.7)}
\def\eqDAh{(4.8)}
\def\thDAca{4.3}
\def\eqDAla{(4.9)}
\def\eqDAlb{(4.10)}
\def\eqDAk{(4.11)}
\def\thDAc{4.4}
\def\eqDAi{(4.12)}
\def\eqDAj{(4.13)}
\def\eqDAm{(4.14)}
\def\eqDAn{(4.15)}
\def\eqDAo{(4.16)}
\def\thDAd{4.5}
\def\eqDApaa{(4.17)}
\def\eqDApab{(4.18)}
\def\eqDAq{(4.19)}
\def\thDBxa{4.6}
\def\thDBa{4.7}
\def\eqDBe{(4.20)}
\def\eqDBf{(4.21)}
\def\eqDBg{(4.22)}
\def\eqDBh{(4.23)}
\def\eqKd{(4.24)}
\def\thDBb{4.8}
\def\eqDBi{(4.25)}
\def\eqDBr{(4.26)}
\def\eqDBia{(4.27)}
\def\eqDBl{(4.28)}
\def\eqDBib{(4.29)}
\def\eqDBk{(4.30)}
\def\eqDBcrx{(4.31)}
\def\eqDBntx{(4.32)}
\def\eqDBj{(4.33)}
\def\thDBd{4.9}
\def\eqKe{(4.34)}
\def\eqKf{(4.35)}
\def\eqDBo{(4.36)}
\def\eqDBp{(4.37)}
\def\eqDBq{(4.38)}
\def\eqLb{(4.39)}
\def\thDBe{4.10}
\def\eqDBa{(4.40)}
\def\thDCf{4.11}
\def\eqDCo{(4.41)}
\def\eqDCoa{(4.42)}
\def\thDCfa{4.12}
\def\eqDCn{(4.43)}
\def\eqDCnz{(4.44)}
\def\eqDBsa{(4.45)}
\def\thDBc{4.13}
\def\eqDBma{(4.46)}
\def\eqDBm{(4.47)}
\def\eqDBn{(4.48)}
\def\eqKg{(4.49)}
\def\eqKh{(4.50)}
\def\eqKi{(4.51)}
\def\eqDCb{(4.52)}
\def\thDCa{4.14}
\def\eqDCaa{(4.53)}
\def\eqDCab{(4.54)}
\def\eqDCac{(4.55)}
\def\eqDCad{(4.56)}
\def\eqDCc{(4.57)}
\def\thDCb{4.15}
\def\eqDCd{(4.58)}
\def\eqDCe{(4.59)}
\def\eqDCda{(4.60)}
\def\eqDCj{(4.61)}
\def\thDCc{4.16}
\def\eqDCk{(4.62)}
\def\eqDCl{(4.63)}
\def\thDCe{4.17}
\def\eqDCm{(4.64)}
\def\thDCd{4.18}
\def\eqDresc{(4.65)}
\def\eqDCf{(4.66)}
\def\eqDCg{(4.67)}
\def\eqDCh{(4.68)}
\def\eqDCha{(4.69)}
\def\eqBBg{(4.70)}
\def\eqDCi{(4.71)}
\def\eqDCia{(4.72)}
\def\eqDCxa{(4.73)}
\def\eqDCxb{(4.74)}
\def\eqDCmc{(4.75)}
\def\thDCg{4.19}
\def\eqDCma{(4.76)}
\def\eqDCmb{(4.77)}
\def\thDCh{4.20}
\def\eqDCna{(4.78)}
\def\eqKj{(4.79)}
\def\eqKk{(4.80)}
\def\eqKl{(4.81)}
\def\eqDCnc{(4.82)}
\def\exDCpfa{(4.83)}
\def\thDCi{4.21}
\def\eqDCya{(4.84)}
\def\eqCDdub{(4.85)}
\def\eqDCyb{(4.86)}
\def\thCDsymm{4.22}
\def\eqCDsymm{(4.87)}
\def\eqGa{(5.1)}
\def\appAxb{(\hbox {A.}1)}
\def\appAxc{(\hbox {A.}2)}
\def\eqDBza{(\hbox {A.}3)}
\def\appAxd{(\hbox {A.}4)}
\def\appAxe{(\hbox {A.}5)}
\def\appAxf{(\hbox {A.}6)}
\def\appAxca{(\hbox {A.}7)}
\def\appAxff{(\hbox {A.}8)}
\def\appAxg{(\hbox {A.}9)}
\def\appAxh{(\hbox {A.}10)}
\def\appAxxv{(\hbox {A.}11)}
\def\appAzx{(\hbox {A.}12)}
\def\appAxz{(\hbox {A.}13)}
\def\appAxpr{(\hbox {A.}14)}
\def\appAxrs{(\hbox {A.}15)}
\def\appAya{(\hbox {A.}16)}
\def\thappAa{\hbox {A.}1}
\def\KPscoj{(\hbox {A.}17)}
\def\appAyc{(\hbox {A.}18)}
\def\appAxin{(\hbox {A.}19)}
\def\appApfl{(\hbox {A.}20)}
\def\appAflp{(\hbox {A.}21)}
\def\eqqqqone{(\hbox {A.}22)}
\def\appAeqqq{(\hbox {A.}23)}
\def\eqDCmf{(\hbox {A.}24)}
\def\appAmo{(\hbox {A.}25)}
\def\eqDCmh{(\hbox {A.}26)}
\def\eqDCmg{(\hbox {A.}27)}
\def\eqAPAaa{(\hbox {B.}1)}
\def\thAPAb{\hbox {B.}1}
\def\eqAPAca{(\hbox {B.}2)}
\def\eqAPAcb{(\hbox {B.}3)}
\def\thAPAa{\hbox {B.}2}
\def\eqAPAza{(\hbox {B.}4)}
\def\eqAPAzb{(\hbox {B.}5)}
\def\eqKm{(\hbox {B.}6)}
\def\eqKn{(\hbox {B.}7)}
\def\thAPAc{\hbox {B.}3}
\def\eqAPAba{(\hbox {B.}8)}
\def\eqAPAaa{(\hbox {B.}9)}
\def\eqAPAab{(\hbox {B.}10)}
\def\eqAPAa{(\hbox {B.}11)}
\def\eqAPAb{(\hbox {B.}12)}
\def\eqAPAc{(\hbox {B.}13)}
\def\eqAPAd{(\hbox {B.}14)}
\def\eqAPAe{(\hbox {B.}15)}
\def\eqAPAf{(\hbox {B.}16)}
\def\eqKo{(\hbox {C.}1)}
\def\KPcmtwo{(\hbox {C.}2)}
\def\KPcmfour{(\hbox {C.}3)}
\def\KPcmone{(\hbox {C.}4)}
\def\appXa{(\hbox {C.}5)}
\def\appXb{(\hbox {C.}6)}
\def\eqKp{(\hbox {C.}7)}
\def\KPCls{(\hbox {C.}8)}
\def\KPcrlem{\hbox {C.}1}
\def\KPCgenla{(\hbox {C.}9)}
\def\KPCgenlb{(\hbox {C.}10)}
\def\KPcmlma{\hbox {C.}2}
\def\KPCprlem{(\hbox {C.}11)}
\def\KPcmlsaa{(\hbox {C.}12)}
\def\KPhtdef{(\hbox {C.}13)}
\def\KPsymth{\hbox {C.}3}
\def\KPCcnsum{(\hbox {C.}14)}
\def\eqKq{(\hbox {C.}15)}
\def\eqKr{(\hbox {C.}16)}
\def\eqKs{(\hbox {C.}17)}
\def\eqKt{(\hbox {C.}18)}
\def\KPcent{\hbox {C.}4}
\def\KPtoshow{(\hbox {C.}19)}
\def\eqKu{(\hbox {C.}20)}
\def\eqKv{(\hbox {C.}21)}
\def\eqKw{(\hbox {C.}22)}
\def\KPcmab{(\hbox {C.}23)}
\def\eqKx{(\hbox {C.}24)}
\def\KPcmgenb{(\hbox {C.}25)}
\def\eqKy{(\hbox {C.}26)}
\def\eqKz{(\hbox {C.}27)}
\def\eqKaa{(\hbox {C.}28)}

\def\virpq{\hbox{{\rm Vir}}_{p,q}}
\def\vir{\hbox{{\rm Vir}}}
\def\tvirq{\widetilde{\hbox{{\rm Vir}}}_{q}}
\def\ni{\noindent}
\def\frac#1#2{{{#1}\over {#2}}}
\def\sqrtN#1{{\!\!\root{N}\of{#1}}}

\def\eql{~=~}
\def\hw{h}
\def\scr{\scriptstyle}
\def\mod{\hbox{\rm mod}\,}
\def\bin#1#2{\left[{#1\atop #2}\right]}
\def\qeii{QEI$\!$I}
\def\NO#1{:\!{#1}\!:}
\def\wT{\widetilde{T}}
\def\wLa{\widetilde{\La}}
\def\wimath{\widetilde{\imath}}
\def\wcH{\widetilde{\cH}}
\def\wM{\widetilde{M}}
\def\qN{q=\sqrtN{1}}
\def\height{{\rm ht}}

\magnification=1200
\hfuzz=20pt
\hsize=6.5truein
\vsize=9.0truein

\nopagenumbers
\pageno=0


\bref\ABF{
G.E.~Andrews, R.J.~Baxter and P.J.~Forrester,
{\it Eight-vertex SOS model and generalized Rogers-Ramanujan-type
identities},
J.\ Stat.\ Phys.\ {\bf 35} (1984) 193-266.}

\bref\AKOSa{
H.~Awata, H.~Kubo, S.~Odake and J.~Shiraishi,
{\it Quantum $\cW_N$ algebras and Macdonald polynomials},
\CMP{179} (1996) 401-416,
{\tt q-alg/9508011}.}

\bref\AKOSc{
H.~Awata, H.~Kubo, S.~Odake and J.~Shiraishi,
{\it Quantum deformation of the $W_N$ algebra}, 
{\tt q-alg/9612001}.}

\bref\AKOSb{
H.~Awata, H.~Kubo, S.~Odake and J.~Shiraishi,
{\it Virasoro-type symmetries in solvable models},
to appear in the CRM series in Mathematical Physics, (Springer Verlag),
{\tt hep-th/9612233}.} 

\bref\BMP{
P.~Bouwknegt, J.~McCarthy and K.~Pilch,
{\it The $\cW_3$ algebra; modules, semi-infinite cohomology and
BV-algebras},
Lect. Notes in Physics Monographs, {\bf m42}
(Springer Verlag, Berlin, 1996).}

\bref\BS{
P.~Bouwknegt and K.~Schoutens,
{\it Spinon decomposition and Yangian structure of $\hsln$ modules}, in
``Geometric Analysis and Lie Theory in Mathematics and Physics'',
Lecture Notes Series of the Australian Mathematical Society,
(Cambridge University Press, Cambridge, 1997),
{\tt q-alg/9703021}.}

\bref\DFJMN{
B.~Davies, O.~Foda, M.~Jimbo, T.~Miwa and A.~Nakayashiki,
{\it Diagonalization of the XXZ Hamiltonian by vertex operators},
\CMP{151} (1993) 89-153.}

\bref\FeFr{
B.~Feigin and E.~Frenkel,
{\it Affine Kac-Moody algebras at the critical level and Gelfand-Dikii
algebras}, \IJMP{A7} (Suppl.\ A1) (1992) 197-215.}

\bref\FF{
B.~Feigin and E.~Frenkel,
{\it Quantum $\cW$-algebras and elliptic algebras},
\CMP{178} (1996) 653--678, 
{\tt q-alg/9508009}.}

\bref\FJMOP{
B.~Feigin, M.~Jimbo, T.~Miwa, A.~Odesskii and Y.~Pugai,
{\it Algebra of screening operators for the deformed $W_n$ algebra},
{\tt q-alg/9702029}.}

\bref\FRa{
E.~Frenkel and N.~Reshetikhin,
{\it Quantum affine algebras and deformations of the Virasoro 
and $\cW$-algebras},
\CMP{178} (1996) 237--264,
{\tt q-alg/9505025}.}

\bref\FRb{
E.~Frenkel and N.~Reshetikhin,
{\it Towards deformed chiral algebras},
{\tt q-alg/9706023}.}

\bref\FRc{
E.~Frenkel and N.~Reshetikhin,
{\it Deformations of $\cal W$-algebras associated to simple Lie algebras},
{\tt q-alg/9708006}.}

\bref\FRS{
E.~Frenkel, N.~Reshetikhin and M.A.~Semenov-Tian-Shansky,
{\it Drinfeld-Sokolov reduction for difference operators and
deformations of $\cW$-algebras, I. The case of Virasoro algebra},
{\tt q-alg/9704011}.}

\bref\Gar{
A.M.~Garsia, {\it Orthogonality of Milne's polynomials
and raising operators}, Discr.\ Math.\ {\bf 99} (1992) 247-264.}

\bref\Gen{
G.~Gentile, 
{\it Osservazioni sopra le statistiche intermedie},
Nuovo Cimento {\bf 17} (1940) 493-497.}

\bref\HYa{
B.-Y.~Hou and W.-L.~Yang,
{\it A $\hbar$-deformed Virasoro algebra as a hidden symmetry of the 
restricted Sine-Gordon model},
{\tt hep-th/9612235}.}

\bref\HYb{
B.-Y.~Hou and W.-L.~Yang,
{\it An $\hbar$-deformation of the $W_N$ algebra and its vertex operators},
J.\ Phys.\ {\bf A}: Math.\ Gen.\ {\bf 30} (1997) 6131-6145,
{\tt hep-th/9701101}.}

\bref\Jia{
N.~Jing, 
{\it Vertex operators, symmetric functions and the spin group $\Ga_n$},
J.~Algebra {\bf 138} (1991) 340-398.}

\bref\Jib{
N.~Jing, 
{\it Vertex operators and Hall-Littlewood symmetric functions},
\AdM{87} (1991) 226-248.}

\bref\JLMP{
M.~Jimbo, M.~Lashkevich, T.~Miwa and Y.~Pugai,
{\it Lukyanov's screening operators for the deformed Virasoro algebra},
Phys.\ Lett.\ {\bf 229A} (1997) 285-292,
{\tt hep-th/9607177}.}

\bref\LLTa{
A.~Lascoux, B.~Leclerc and J.-Y.~Thibon, {\it Fonctions de 
Hall-Littlewood et polyn\^omes de Kostka-Foulkes aux racines de l'unit\'e},
C.R.\ Acad.\ Sci.\ Paris {\bf 316} (1993) 1-6.}

\bref\LLTb{
A.~Lascoux, B.~Leclerc and J.-Y.~Thibon, {\it Green polynomials and
Hall-Littlewood functions at roots of unity}, Europ.\ J.\ Combinatorics
{\bf 15} (1994) 173-180.}

\bref\Lit{
D.E.~Littlewood,
{\it On certain symmetric functions},
Proc.\ London Math.\ Soc.\ {\bf 11} (1961) 485-498.}

\bref\LPa{
S.~Lukyanov and Y.~Pugai,
{\it Bosonization of ZF algebras: Direction toward deformed Virasoro
algebra}, J.\ Exp.\ Theor.\ Phys.\ {\bf 82 } (1996) 1021-1045,
{\tt hep-th/9412128}.}

\bref\LPb{
S.~Lukyanov and Y.~Pugai,
{\it Multi-point local height probabilities in the integrable RSOS
model}, \NPB{473} (1996) 631-658,
{\tt hep-th/9602074}.}

\bref\MD{
I.G.~Macdonald,
{\it Symmetric functions and Hall polynomials},
(Oxford University Press, Oxford, 1995).}

\bref\Mila{
S.C.~Milne, {\it Classical partition functions and
the $U(n+1)$ Rogers-Selberg identity}, Discr.\ Math.\ {\bf 99} (1992) 199-246.}

\bref\Milb{
S.C.~Milne, {\it The $C_\ell$ Rogers-Selberg identity},
SIAM J.\ Math.\ Anal.\ {\bf 25} (1994) 571-595.}

\bref\Mob{
A.O.~Morris, 
{\it On an algebra of symmetric functions},
Quart.\ J.\ Math.\ Oxford {\bf 16} (1965) 53-64.}

\bref\Pol{
A.~Polychronakos, {\it Path integrals and parastatistics}, 
\NPB{474} (1996) 529-539, {\tt hep-th/9603179}.}

\bref\Sch{
K.~Schoutens, {\it Exclusion statistics in conformal field theory
spectra}, Phys.\ Rev.\ Lett.\ {\bf 79} (1997) 2608-2611,
{\tt cond-mat/9706166}.}

\bref\SS{
A.M.~Semenov-Tian-Shansky and A.V.~Sevostyanov,
{\it Drinfeld-Sokolov reduction for difference operators and
deformations of $\cW$-algebras, II. General semisimple case},
{\tt q-alg/9702016}.}

\bref\SKAO{
J.~Shiraishi, H.~Kubo, H.~Awata and S.~Odake,
{\it A quantum deformation of the Virasoro algebra and the Macdonald
symmetric functions},
\LMP{38} (1996) 33--51,
{\tt q-alg/9507034}.}

%
%
\rightline{\tt\date}
\vskip2cm
\centerline{\bf THE DEFORMED VIRASORO ALGEBRA AT ROOTS OF UNITY}\bigskip
\vskip1cm

\centerline{Peter BOUWKNEGT$\,^{1}$ and Krzysztof PILCH$\,^{2}$}
\bigskip

\centerline{\sl $^1$ Department of Physics and Mathematical Physics}
\centerline{\sl University of Adelaide}
\centerline{\sl Adelaide, SA~5005, AUSTRALIA}
\bigskip

\centerline{\sl $^2$ Department of Physics and Astronomy}
\centerline{\sl University of Southern California}
\centerline{\sl Los Angeles, CA~90089-0484, USA}
\medskip
\vskip1.5cm

\centerline{\bf ABSTRACT}\bigskip

\vbox{\leftskip 2.0truecm \rightskip 2.0truecm
\noindent 
We discuss some aspects of the representation theory of the 
deformed Virasoro algebra $\virpq$.  In particular, we give a proof of
the formula for the 
Kac determinant and then determine the center of $\virpq$ for
$q$ a primitive $N$-th root of unity.  We derive explicit expressions
for the generators of the center in the limit $t=qp^{-1}\rightarrow
\infty$ and elucidate the connection to the Hall-Littlewood 
symmetric functions.  Furthermore, we argue that for
$q=\sqrtN{1}$ the algebra describes `Gentile statistics' of order
$N-1$, i.e., a situation in which at most $N-1$ particles can occupy
the same state.}

\vfil
\leftline{ADP-97-22/M53}
\leftline{USC-97/17}
\line{{\tt q-alg/9710026}\hfil October 1997}

\eject


\baselineskip=1.5\baselineskip

\footline{\hss \tenrm -- \folio\ -- \hss}

\newsec{Introduction}
 
In recent years it has been realized (see, in particular, [\DFJMN])
that the theory of off-critical integrable models of statistical
mechanics is intimately connected to the theory of
infinite-dimensional quantum algebras and that such theories can be
studied in close parallel to their critical counterparts, i.e.,
conformal field theories.

Algebras of particular interest are the so-called deformed Virasoro
algebra, $\virpq$, introduced in [\FRa,\SKAO] (see [\AKOSb] for a 
review), their higher rank 
generalization, the deformed $\cW$-algebras [\FF,\AKOSa,\AKOSc,\FRc], 
as well as their linearized
versions [\HYa,\HYb].  In particular, it has been argued [\LPa,\LPb] that
the deformed Virasoro algebra plays the role of the dynamical symmetry
algebra in the Andrews-Baxter-Forrester RSOS models [\ABF].

For generic values of the deformation parameters\foot{Here, `generic
values of the deformation parameters' will always stand for `not a root
of unity.'} the representation theory of the deformed Virasoro 
algebra $\virpq$ closely parallels that of the undeformed Virasoro
algebra, as manifested, e.g., by the Kac determinant formula [\SKAO],
Drinfel'd-Sokolov reductions [\FRS,\SS] and the existence of Felder type
free field resolutions [\LPb,\JLMP,\FJMOP].

In this paper we will discuss some aspects of the representation
theory of $\virpq$.  First we complete the proof of the Kac determinant 
conjectured in [\SKAO], where also the hardest part of the proof,
the construction of a sufficient number of vanishing lines, was already
established.  Then we proceed to discuss the 
representation theory for $q$ a primitive $N$-th root of unity, i.e., a
complex number $q$ such that $q^N=1$ and $q^k\neq1$ for all $0<k<N$.
(Henceforth, we use the notation $q=\sqrtN{1}$ for $q$ a primitive
$N$-th root of unity.)  As was already observed in [\SKAO], it follows
from the Kac determinant formula that, for $q=\sqrtN{1}$, Verma 
modules contain
many singular vectors irrespective of the highest weight.
We first analyze the case $q=-1$ in detail and find a close relation 
of $\virpq$ to the free fermion algebra.
For $q=\sqrtN{1}$, $N>2$, we give a generating series for all singular 
vectors.  We derive explicit expressions of all singular vectors 
(for generic highest weight) in the limit $t=q/p\to \infty$. 
For generic $t$, the
singular vectors are deformations of these as we illustrate in
various examples.
We show that the existence of these generic singular vectors is a 
consequence of the fact that for $q=\sqrtN{1}$ (and $t\to\infty$)
the algebra $\virpq$ has a large center.  We compute this center by
exploiting the isometry from the Verma module to the Hall-Littlewood
symmetric polynomials [\Jib].

The paper is organized as follows.  In Section 2 we recall the
definition of the deformed Virasoro algebra $\virpq$ and its 
highest weight modules,
prove some simple properties
and discuss the free field realization.
In Section 3 we prove the Kac determinant formula and
in Section 4 we discuss the representation
theory of $\virpq$ for $q=\sqrtN{1}$.  This section is divided
into three parts.  First we discuss the case of $q=-1$
and then proceed to the general case of $q=\sqrtN{1}$ for all $N\in\NN$.
Finally, we make the previous analysis more explicit in
the limit $t=q/p\to \infty$.  We conclude with some general comments in
Section 5.  In particular, we
point out an interesting
relation between $\virpq$ at $q=\sqrtN{1}$, and so-called Gentile
statistics of order $N-1$.  The latter is 
a generalization of Fermi statistics
($N=2$) in which at most $N-1$ particles can occupy the same state [\Gen].

Three appendices follow.  In Appendix A we give basic definitions and
summarize some results from the theory of symmetric functions that are
used throughout the paper.
In Appendix B we provide some explicit 
examples of singular vectors at $q=\sqrtN{1}$ for $N=3,4$, and in Appendix
C we establish some elementary identities for the products of
generators of $\virpq$ in the limit $t\to\infty$ which, when specialized
to $q=\sqrtN{1}$, yield another derivation of the center of $\virpq$.

\newsec{The deformed Virasoro algebra $\virpq$}

\subsec{Definition}

The deformed Virasoro algebra [\SKAO,\AKOSb], $\virpq$,  is defined 
to be the associative algebra generated by $\{T_n,\ n\in\ZZ\}$, 
with relations
\eqn\eqBAa{
\sum_{l\geq0} f_l \left(T_{m-l}T_{n+l} - T_{n-l}T_{m+l} \right)
  \eql c_m \de_{m+n,0}\,,
}
where $f_l$ is determined through
\eqn\eqBAb{\eqalign{
f(x) \equiv \sum_{l\geq0} f_l x^l 
& \eql \exp\Biggl( \sum_{n\geq 1} { (1-q^n)(1-t^{-n})\over (1+p^n)}
  {x^n\over n} \Biggr) \cr
& \eql {1\over 1-x} { (qx;p^2)_\infty (q^{-1}px;p^2)_\infty 
\over (qpx;p^2)_\infty (q^{-1}p^2x;p^2)_\infty} \,,\cr}
}
and
\eqn\eqBAc{
c_m \eql \ze (p^m - p^{-m})\,,\qquad 
\ze \eql -{(1-q)(1-t^{-1})\over 1-p}\,.
}
Here, $p,q\in\CC$ with $p$ not a root of $-1$. The series 
\eqBAb\ are to be understood as formal power series and the equality
holds in the region where both converge. 
For convienience we have introduced a third parameter $t$ by $t=qp^{-1}$.
Also
\eqn\eqLa{
(x;q)_M \eql \prod_{k=1}^M\ (1-x q^{k-1})\,,\qquad (q)_M \eql (q;q)_M\,.
}

In terms of formal series
\eqn\eqBAd{
T(z) \eql \sum_{n\in\ZZ} T_n z^{-n}, \qquad \de(z) = \sum_{n\in\ZZ} z^n\,,
}
the relations \eqBAa\ read
\eqn\eqBAe{
f({w\over z}) T(z) T(w) - f({z\over w})T(w)T(z) \eql
  \ze \left( \de(p  {w\over z}) - 
  \de(p^{-1}  {w\over z}) \right)\,,
}

For future convenience we also introduce a generating series for 
$c_m = -c_{-m}$
\eqn\eqBAf{
c(x) ~\equiv~ \sum_{m\geq0} c_m x^m \eql \ze 
\left( {1\over 1-px} - {1\over 1-p^{-1}x}\right)\,.
}

Note that the algebra $\virpq$ is invariant under 
$(p,q)\leftrightarrow (p^{-1},qp^{-1})$ and $(p,q)\leftrightarrow 
(p^{-1},q^{-1})$ and carries a $\CC$-linear anti-involution $\om$ defined
by 
\eqn\eqBAfa{
\om (T_{n}) \eql T_{-n}\,.
}

{}From the second expression for $f(x)$ in \eqBAb\ one can easily derive
the recurrence relation
\eqn\eqBAg{
f(x) f(px) \eql { (1-qx) (1-q^{-1}px) \over (1-x) (1-px)  }
 \eql 1 + \ze pq^{-1} \left( {1\over 1-x} - {1\over 1-px} \right)\,,
}
which in turn uniquely determines $f(x)$.

\ni {\bf Remark.} We have defined 
$\virpq$ for any $p\in\CC$ by means of \eqBAa\ -- \eqBAc\
considered as formal power series, as long as $p$ is not a root of
$-1$.  For $p$ a root of $-1$, the expressions \eqBAb\ are
ill-defined. In this case one may still define $\virpq$ in terms of
a solution of the recurrence relation
\eqBAg.  However, such a solution is not unique as can be seen by 
taking different limits in \eqBAb\ to obtain inequivalent solutions to
\eqBAg. It is rather straightforward to show by iterating \eqBAg\
that, for $p$ an $N$-th root of $-1$, a necessary and sufficient
condition for the existence of a solution is $q^{2N}=1$.  We
have not found a succinct way to classify all the solutions that arise
then. In this paper we will focus on the case where $p$ is not a root of
$-1$, except for brief remarks in Sections 2.2 and 4.1.\medskip

The algebra $\virpq$ can be considered to be a deformation of
the Virasoro algebra $\vir$
\eqn\eqBAh{
[L_m,L_n] \eql (m-n) L_{m+n} + {c\over 12} m(m^2-1) \de_{m+n,0}\,,
}
in the following sense:
\thm\thBAa
\proclaim Theorem \thBAa\ [\SKAO]. In the limit $q = e^\hbar \to1$ 
with $t=q^\be$ ($\be$ fixed) we have 
\eqn\eqBAi{
T(z) \eql 2 + \be \left( L(z) + {(1-\be)^2 \over 4\be} \right)\hbar^2 + 
  \cO(\hbar^4)\,,
}
where $L(z) = \sum_{n\in\ZZ} L_n z^{-n}$ satisfies the Virasoro algebra 
\eqBAh\ with
\eqn\eqBAj{
c \eql 1 - {6(1-\be)^2\over \be}\,.
}

\subsec{Modules}

In this paper we will consider the class of modules of $\virpq$ that
are analogues, e.g., deformations, of the highest weight
modules of the undeformed algebra. In particular, we are interested in
studying those modules by means of standard techniques based on
characters and contravariant forms. For that reason we extend $\virpq$
by a derivation $d$ satisfying
\eqn\eqCaa{
[d,T_n] \eql -n T_n\,,
}
and define the category $\cO$ of $\virpq$ modules as the set of
$d$-diagonalizable modules $V=\amalg_{z\in\CC} V_{(z)}$, 
such that each $d$-eigenspace $V_{(z)}$ is finite dimensional.  Let
$P(V) = \{z\in\CC\,|\, V_{(z)}\neq0\}$.  We also require that there exists
a finite set $z_1,\ldots,z_s\in\CC$ such that 
$P(V) \subset \cup_{i=1}^s \{ z_i + \NN\}$.
Then, the character of a $\virpq$ module $V\in\cO$ 
is defined by
\eqn\eqCe{
{\rm ch}_V(x) \eql {\rm Tr}\ x^d \eql \sum_{z\in\CC} {\rm dim}(V_{(z)})\,
x^z\,.
}
One can show that so defined category $\cO$ contains, in particular,
highest weight modules, that are defined in the usual manner and
include, among others, Verma modules and their (irreducible) quotients.

The Verma module $M(h)$ [\SKAO,\FF], with 
highest weight state $|h\rangle$ satisfying $T_0|\hw\rangle
=\hw |\hw\rangle$ and
$T_n|h\rangle=0$, $n>0$, has a basis indexed by 
partitions $\la=(\la_1,\la_2,\ldots)$, $\la_1\geq\la_2\geq \ldots >0$,
i.e., a basis of $M(\hw)_{(n)}$ is
given by the vectors 
\eqn\eqCAf{
|\la;\hw\rangle ~\equiv~
T_{-\la} |\hw\rangle ~\equiv~ T_{-\la_1} \ldots T_{-\la_\ell} |\hw\rangle
\eql T_{-n}^{m_n(\la)} \ldots T_{-1}^{m_1(\la)} |\hw\rangle\,,
}
where $\la$ runs through all partitions of $n$.  We use the notation
$\la\vdash n$.  
Furthermore, $m_i(\la) = \#\{j \, : \,  \la_j = i\}$ 
denotes the number of parts of length $i$ in $\la$.  For $\la\vdash n$
we also use $|\la|=n$ and define the length of $\la$ by $\ell(\la) =
\sum_{i\geq1} m_i(\la)$.
The following two orderings are used; the (reverse)
lexicographic ordering, i.e.,
$\la > \mu$ if the 
first non-vanishing difference $\la_1-\mu_1$, $\la_2-\mu_2$,
$\ldots$, is positive, and the natural (partial) ordering, i.e.,
$\la\succeq \mu$ if $\sum_{j=1}^i \la_j \geq
\sum_{j=1}^i \mu_j$ for all $i\geq1$.

{}From the previous discussion it follows that the character of the 
Verma module $M(\hw)$ is given by
\eqn\eqCAfa{
{\rm ch}_M(x) \eql \prod_{n\geq1}\left( {1\over 1-x^n}\right) \eql
\sum_{n\geq0}\ p(n) \,x^n\,,
}
where $p(n)$ is the number of partitions of $n$.

The anti-involution $\om$ of $\virpq$ (see \eqBAfa) determines a
bilinear contravariant form (`Shapovalov form') on $M(\hw)$,
uniquely defined by
\eqn\eqCAg{
G_{\la\mu} ~\equiv~ \vev{ \la;\hw | \mu;\hw } \eql 
  \vev{ \hw | \om(T_{-\la}) T_{-\mu} | \hw}\,,
}
and 
\eqn\eqCAh{
\vev{ \hw | \hw } \eql 1\,.
}
In particular we have $G_{\la\mu}=0$ for $|\la|\neq |\mu|$.  
In Section 3 we compute the so-called Kac determinant $G^{(n)}$, i.e.,
the determinant of the form $G_{\la\mu}$ on $M(\hw)_{(n)}$.

As an aside, one
might wonder whether the category $\cO$ contains any finite dimensional
irreducible representations of $\virpq$.  Clearly, $\virpq$ is
abelian for $q=1$ and/or $t=1$ and therefore has a wealth of one dimensional
irreducible representations in these cases.  The complete result is: 
\thm\thCc
\proclaim Theorem \thCc. The only irreducible finite dimensional 
$\virpq$ modules in $\cO$ are one dimensional and occur for 
\item{(i)} $q=1$ and/or $t=1$, and $h\in\CC$ arbitrary.
\item{(ii)} $p=q^{{1\over3}}$ or $p=q^{-{1\over2}}$ and
$h^2 = p^{-1} (1+p)^2$.
\medskip

\ni
{\bf Remark.} In case (i) we have $c(x)=0$ and $f(x)=1$.  Case (ii)
corresponds to a $c(x)\not=0$ deformation of the $c=0$ Virasoro module
with $h^2 = h^2_{1,1}$ (cf.\ \eqBAj\ and (3.6)).
\medskip

\proof Let $|h\rangle$ be the highest weight
state of a finite dimensional irreducible module in $\cal O$. Define
$a_n=\langle h|T_nT_{-n}|h\rangle$, $n\geq 0$. Then $a_0=h^2$, and, as
follows  from the commutation relations
\eqBAa,
\eqn\KPbasiden{
a_n+\sum_{l=1}^na_{n-l}f_l \eql c_n\,,\qquad n\geq 1\,.
}
Since the module is finite dimensional, we must have $a_n=0$ for
$n$ sufficiently large, say $n>n_0$. Multiplying 
\KPbasiden\ by $x^n$ and summing over $n\geq n_0+1$, we find 
\eqn\KPfinid{
\big(\sum_{i=0}^{n_0} a_ix^i\,\big)f(x) \eql c(x)+a_0\,.
}
In arriving at \KPfinid\ we have also used \KPbasiden\ with $n\leq
n_0$ to simplify some intermediate expressions.

First consider the case $c(x)=0$. This can occur only for the
following values of the deformation parameters:  $q=1$, $p$
arbitrary;  $p=q$, $q$ arbitrary;  $p=-1$, $q$ arbitrary. (In
view of the remark in Section~2.1, the last case is in fact covered by
the first two.) We conclude that for $c(x)=0$, $\virpq$ is an abelian
algebra and its irreducible modules are one dimensional as described
by case (i).

Now, suppose $c(x)\not=0$. In this case not all $a_n$, $n\geq 0$,
vanish and $f(x)$ is a rational function. It follows from the
recurrence relation \eqBAg\ that
$\lim_{x\rightarrow\infty}f(x)=\pm1$. Since $\lim_{x\rightarrow
\infty} c(x)=0$, we must have $a_n=0$ for $n\geq 1$, so that
\KPfinid\ becomes 
\eqn\KPfandc{
h^2(f(x)-1)\eql c(x)\,.
}
In terms of modes this is simply $h^2f_n=c_n$.  We also deduce from
$a_1=0$ that
\eqn\eqKa{
h^2 \eql{c_1\over f_1}\eql {(1+p)^2\over p}\,.
}
After solving  \KPfandc\ for $f(x)$, we find that 
\eqBAg\ holds if and only if $q$ and $p$ 
satisfy one of the conditions in (i) or
(ii). In particular, in case (ii),
\eqn\eqBAxa{
f(x)\eql {(1-p^{-2}x)(1-p^2x)\over (1-p^{-1}x)(1-px)}\,.}

We must still show that the resulting module is one dimensional.  Let
$m_0$ be the smallest $m>0$ for which $T_{-m}|h\rangle\not=0$. We have
shown already that $T_{m_0}T_{-m_0}|h\rangle=0$. For
$n=1,\ldots,m_0-1$ we find
\eqn\KPsing{
T_nT_{-m_0}|h\rangle
\eql -\sum_{l=1}^{m_0-1}T_{n-l}T_{-m_0+l}|h\rangle-h f_{m_0}
T_{n-m_0}|h\rangle \eql 0\,.
}
Thus $T_{-m_0}|h\rangle$ is a singular vector and must vanish. \Box\medskip


\subsec{The $q$-Heisenberg algebra $\cH_{p,q}$ and the 
free field realization of $\virpq$}

The $q$-Heisenberg algebra, $\cH_{p,q}$, is the associative algebra 
with generators $\{\al_n,\ n\in\ZZ\}$ and relations 
\eqn\eqBBa{
[\al_m,\al_n] \eql m\, { (1-q^m)(1-t^{-m})\over 1+p^m }\, \de_{m+n,0}\,.
}
Let $U(\cH_{p,q})_{\rm loc}$ denote the local completion of $\cH_{p,q}$
(see [\FeFr]).  Furthermore, let $\om_{\cH}$ denote the $\CC$-linear 
anti-involution of $U(\cH_{p,q})_{\rm loc}$ defined
by 
\eqn\eqBBb{ \eqalign{
\om_{\cH}(\al_m) & \eql p^{-m} \al_{-m}\,, \qquad m\neq0\,,\cr
\om_{\cH}(q^{\al_0}) & \eql p q^{-\al_0}\,.\cr}
}
Again, we can add a derivation $d$ to $\cH_{p,q}$ defined by
\eqn\eqBBaa{
[d,\al_m] \eql -m \al_m\,,\qquad m\in\ZZ\,.
}
For any $\al\in\CC$ we have an $\cH_{p,q}$ module $F(\al)$, the 
so-called Fock space, which is irreducible for generic $p$ and $q$, and 
decomposes as 
\eqn\eqBBab{
F(\al) \eql \amalg_{n\geq0} \ F(\al)_{(n)}\,,
}
under the action of $d$.  The module $F(\al)_{(n)}$ has a basis indexed 
by partitions $\la\vdash n$
\eqn\eqBBac{
|\la,\al\rangle ~\equiv~ \al_{-\la} |\al\rangle 
 ~\equiv~ \al_{-\la_1}\ldots \al_{-\la_\ell} 
 |\al\rangle\,,
}
and where the highest weight vector (the `vacuum') satisfies
\eqn\eqBBad{ \eqalign{
\al_0 |\al\rangle & \eql \al |\al\rangle \,,\cr
\al_m |\al\rangle & \eql 0\,,\qquad  m>0\,.\cr}
}
The anti-involution $\om_\cH$ of \eqBBb\ induces a unique contravariant
bilinear form $\vev{-|-}_F$ on $F(\al') \times F(\al)$ such that
$\vev{\al'|\al}=1$, where $\al'$ is determined by $q^{\al'} = 
p q^{-\al}$, i.e.,
\eqn\eqBBae{
g_{\la\mu} ~\equiv~ \vev{\la;\al'|\mu;\al}_F \eql
 \vev{\al' | \om_\cH(\al_{-\la}) \al_{-\mu} | \al}_F\,.
}
We compute $g_{\la\mu}$ explicitly in Section 3.

We now recall the free field realization of $\virpq$.
\thm\thBBa
\proclaim Theorem \thBBa\ [\FF,\SKAO]. We have a homomorphism of 
algebras 
$\imath: \virpq \to U(\cH_{p,q})_{\rm loc}$ defined by
\eqn\eqBBc{
\imath( T(z)) \eql \La^+(z) + \La^-(z)\,,
}
where
\eqn\eqBBd{ \eqalign{
\La^+(z) & \eql p^{-{1\over2}} q^{\al_0} \exp\Biggl(
\sum_{n\geq1} {\al_{-n}\over n} z^n \Biggr) \exp\Biggl( -\sum_{n\geq1}
{\al_n\over n} z^{-n} \Biggr)\,,\cr
\La^-(z) & \eql p^{{1\over2}} q^{-\al_0} \exp\Biggl( -
\sum_{n\geq1} {\al_{-n}\over n} (p^{-1}z)^n \Biggr) \exp\Biggl( \sum_{n\geq1}
{\al_n\over n} (p^{-1}z)^{-n} \Biggr)\,.\cr}
}
Note that we can write
\eqn\eqBBdy{
\La^+(z) \eql \NO{ \La(z) }\,,\qquad \La^-(z) \eql \NO{ \La(p^{-1}z)^{-1} }\,,
}
where 
\eqn\eqBBdz{
\La(z) \eql p^{-{1\over2}} q^{\al_0} \exp \Biggl( - 
\sum_{n\neq0} {\al_n\over n}
z^{-n} \Biggr) \,.
}
Furthermore
\eqn\eqBBe{
\imath \circ \om \eql \om_\cH \circ \imath\,.
}
\medskip

We sketch the proof since some intermediate results will be
useful in later sections.\medskip

\proof By means of standard free field techniques we find
for $|z_1|\gg |z_2|$ and $\ep_i\in\{\pm\}$
\eqn\eqBBf{ 
\La^{\ep_1}(z_1) \La^{\ep_2}(z_2)  \eql f^{\ep_1\ep_2}({z_2\over z_1}) 
\NO{\La^{\ep_1}(z_1) \La^{\ep_2}(z_2)}\,,
}
with
\eqn\eqBBg{
f^{++}(x) \eql f^{--}(x) \eql f(x)^{-1}\,,\qquad f^{+-}(x)\eql f(p^{-1}x)\,,
\qquad f^{-+}(x)\eql f(px)\,,
}
where $f(x)$ is defined in \eqBAb.  Thus we have 
\eqn\eqBBh{
f( {z_2\over z_1} ) T(z_1)T(z_2) \eql \sum_{\ep_i} 
F^{\ep_1\ep_2}({z_2\over z_1}) \NO{\La^{\ep_1}(z_1) \La^{\ep_2}(z_2)}\,,
}
with
\eqn\eqBBi{
F^{++}(x) \eql F^{--}(x) \eql 1\,,\qquad F^{+-}(x)\eql f(x)f(p^{-1}x)\,,
\qquad F^{-+}(x)\eql f(x)f(px)\,.
}
Subtracting the term with $z_1\leftrightarrow z_2$ gives
\eqn\eqBBj{ \eqalign{
f( {z_2\over z_1} ) T(z_1)T(z_2) - f( {z_1\over z_2} ) T(z_2)T(z_1)
\eql &  (F^{+-}(x) - F^{-+}({1\over x})) \NO{ \La^+(z_1)\La^-(z_2)} \cr & +
     (F^{-+}(x) - F^{+-}({1\over x})) \NO{ \La^-(z_1)\La^+(z_2)}\,.\cr}
}
Using \eqBAg\ we find
\eqn\eqBBk{ \eqalign{
F^{+-}(x) - F^{-+}({1\over x}) & \eql \ze ( \de(p^{-1}x) - \de(x) )\,,\cr
F^{-+}(x) - F^{+-}({1\over x}) & \eql \ze ( \de(x) - \de(px)  )\,,\cr}
}
while \eqBBd\ gives
\eqn\eqBBl{
\NO{ \La^+(z)\La^-(pz)} \eql 1\,.
}
This completes the proof.  \Box\medskip

The free field realization $\imath:\virpq \to U(\cH_{p,q})_{\rm loc}$
equips the $U(\cH_{p,q})_{\rm loc}$ module $F(\al)$ with the structure
of a $\virpq$ module.  In fact, from Theorem \thBBa\ it follows
\thm\thBBb
\proclaim Corollary \thBBb.  Let $\al\in\CC$ and $p,q\in\CC$ arbitrary.
Define 
\eqn\eqBBba{
\hw(\al) \eql p^{-{1\over2}} q^{\al} + p^{{1\over2}} q^{-\al}\,.
}
There exists a unique homomorphism $\imath$ of $\virpq$ modules,
$\imath: M(\hw(\al)) \rightarrow F(\al)$, 
such that $\imath( |\hw(\al)\rangle ) = |\al\rangle$.
The homomorphism $\imath$ is an isometry with respect 
to the bilinear contravariant forms defined on $M(\hw)$ and $F(\al)$.
\medskip

\ni {\bf Remark.} Using  \eqBBc\ and \eqBBf\ one can show that
$\imath(T(z))$ satisfies, in the sense of meromorphically
continued products of operators, the following exchange relation
[\FF],
\eqn\eqBBexch{ 
\imath(T(z))\imath(T(w)) \eql S_{TT}({w\over
z})\imath(T(w))\imath(T(z))\,,
}  
where
\eqn\eqBBstt{ 
S_{TT}(x)\eql f({1\over x})f(x)^{-1}\,,
} 
is a solution to the Yang-Baxter equation. This observation is further
developed in [\FRb], where a proposal is made for the definition of a
deformed chiral algebra (DCA). In this formalism the deformed Virasoro
algebra, $\virpq$, is naturally defined as a subalgebra of the DCA
corresponding to the $q$-deformed Heisenberg algebra,
$\cH_{p,q}$. While some of our discussion can be naturally recast in
the language of DCAs, for most of our purposes the algebraic setup of
Section 2.1 is sufficient.


\newsec{The Kac determinant}

An explicit formula for the Kac determinant of the 
Verma modules of $\virpq$ 
was conjectured in [\SKAO].  In this section we present its
proof by using the isometry $\imath: M(\hw(\al)) \to F(\al)$.

Let us first consider the determinant of the bilinear form
on $F(\al)$.  For partitions $\la,\mu\vdash n$ we have
\eqn\eqCBaa{
g^{(n)}_{\la\mu} \eql \vev{\al | \om_{\cH}(\al_{-\la}) \al_{-\mu} }_F
\eql \de_{\la\mu}\ z_\la\ p^{|\la|} \prod_{i=1}^{\ell(\la)} 
{ (1-q^{\la_i})(1-t^{-\la_i})\over 1+p^{\la_i} }\,,
}
where
\eqn\eqCBab{
z_\la \eql \prod_{i\geq1}\  i^{m_i(\la)}\, m_i(\la)!\,.
}
\thm\thCBaa
\proclaim Theorem \thCBaa. 
\eqn\eqCBac{
g^{(n)} ~\equiv~  {\rm det}\ g_{\la\mu} \eql
C_n \prod_{{\scr r,s\geq1}\atop {\scr rs\leq n} }
  \left(p^r {(1-q^r)(1-t^{-r}) \over 1 + p^r }
  \right)^{p(n-rs)} \,,
}
where $C_n = \prod_{\la\vdash n} z_\la$ is a constant independent of 
$p$, $q$ and $q^\al$.
\medskip

In deriving Theorem \thCBaa\ 
we have used the following elementary result.
\thm\thCBab
\proclaim Lemma \thCBab.  Let $f$ be a (complex valued) function on $\NN$.
Then, for all $n\geq0$,
\eqn\eqCBad{
\prod_{\la\vdash n}\ \prod_{i=1}^{\ell(\la)}\ f(\la_i) \eql 
\prod_{\la\vdash n}\ \prod_{i=1}^n\ 
\prod_{j_i=1}^{m_i(\la)}\  f(j_i) \eql
\prod_{{\scr r,s\geq1}\atop {\scr rs\leq n}}\  f(r)^{p(n-rs)}\,.
}
\medskip

\proof Consider the left hand side of \eqCBad. Fix $r\geq1$.  Consider 
all partitions $\la\vdash n$ with a row of length $r$.  Since there are 
$p(n-r)$ partitions with at least one row of length $r$, the first such row 
contributes a factor of $f(r)^{p(n-r)}$.  There are $p(n-2r)$ partitions
with at least two rows of length $r$, so the second row of length $r$
contributes $f(r)^{p(n-2r)}$.  Iterating this we conclude that the
rows of length $r$ in all partitions of $n$ contribute a factor
$\prod_s f(r)^{p(n-rs)}$.  This proves equality of the left hand side
to the right hand side.
The equality of the middle formula to the right hand side
is proved similarly, but now considering rows of length $s$.
\Box\medskip

We are now ready for the main result of this section.
\thm\thCBa
\proclaim Theorem \thCBa.
The Kac determinant of $M(\hw)_{(n)}$ is given by 
\eqn\eqCBa{
G^{(n)} \eql C_n\prod_{{\scr r,s\geq1}\atop {\scr rs\leq n} }
  (\hw^2 - \hw^2_{r,s})^{p(n-rs)} \left( {(1-q^r)(1-t^{-r}) \over 1+p^r }
  \right)^{p(n-rs)} \,,
}
where 
\eqn\eqCBb{
\hw_{r,s} \eql t^{\scr r\over\scr 2}  q^{\scr s\over\scr 2}+
t^{-{\scr r\over\scr 2}} q^{{\scr s\over\scr 2}}
\eql p^{-{\scr r\over\scr 2}} q^{{\scr r-s\over\scr 2}}+
p^{{\scr r\over\scr 2}} q^{\scr  s-r\over\scr 2}\,,
}
and $C_n$ is a constant independent of $p$, $q$ and $\hw$.
\medskip 

\proof  Fix $n\in\NN$.  For partitions $\la,\mu\vdash n$, define 
a matrix $\Pi^{(n)}_{\la\mu}=\Pi^{(n)}_{\la\mu}(p,q,q^\al)$ by
\eqn\eqCBc{
\imath(T_{-\la}) |\al\rangle \eql \sum_{\mu\vdash n} \ 
\Pi^{(n)}_{\la\mu}\ \al_{-\mu}^{\phantom{(}} |\al\rangle \,.
}
and let 
\eqn\eqCBd{
\Pi^{(n)} \eql {\rm det}\, \Pi^{(n)}_{\la\mu}\,.
}
Clearly, $\Pi^{(n)}(p,q,q^\al)$ is a Laurent polynomial in $p$, $q$ and 
$q^\al$.
Using \eqBBb, we have
\eqn\eqCBe{
G^{(n)} \eql \Pi^{(n)}(p,q,pq^{-\al})\ g^{(n)} \ \Pi^{(n)}(p,q,q^\al) \,,
}
where $g^{(n)}$ is given in Theorem \thCBaa.
Now, the crucial step is that, from the explicit construction 
of singular vectors 
in terms of Macdonald polynomials [\SKAO], it is known that $G^{(n)}$ 
has vanishing lines 
\eqn\eqCBf{ \eqalign{
h^2-h_{r,s}^2 \eql &
\left( p^{{ {\scr r-1\over\scr 2}}} q^{{ {\scr s-r\over\scr 2}}} 
q^\al -
  p^{-{ {\scr r-1\over\scr 2}}} q^{{ {\scr r-s\over\scr 2}}} 
q^{-\al} \right)  \cr
& \times \left( p^{-{ {\scr r+1\over\scr 2}}} 
q^{{ {\scr r-s\over\scr 2}}} q^\al -
  p^{{ {\scr r+1\over\scr 2}}} 
q^{{ {\scr s-r\over\scr 2}}} q^{-\al}   \right)\,,\cr}
}
i.e., it has a factor
$$
\prod_{{\scr r,s\geq1}\atop {\scr rs\leq n}} \ 
  (\hw^2 - \hw^2_{r,s})^{p(n-rs)} \,.
$$
Thus, $G^{(n)}$ is given by \eqCBa, up to a Laurent polynomial $C_n$ in 
$p$, $q$ and $q^\al$.  The proof will now be completed if we can show that
$C_n$ is actually a constant.  As a Laurent polynomial in $q^\al$, 
the leading term of
$\Pi^{(n)}(p,q,q^\al)$ is easily computed.  It arises from 
$\Pi^{(n)+}(p,q,q^\al)=
{\rm det}\, \Pi^{(n)+}_{\la\mu}(p,q,q^\al)$,
where 
\eqn\eqCBg{
\La^+_{-\la} |\al\rangle \eql \sum_{\mu} \ 
\Pi^{(n)+}_{\la\mu}\ \al_{-\mu}^{\phantom{(}} |\al\rangle \,.
}
This determinant can be computed by first going to a basis of 
$F(\al)_{(n)}$ given by the vectors
\eqn\eqCBga{
A_{-\la} |\al\rangle ~\equiv~ A_{-\la_1}\ldots A_{-\la_l}|\al\rangle\,,
}
where $\la$ runs over all partitions of $n$ and $A_{-m}$ is defined through
\eqn\eqCBgb{
\exp\Biggl( \sum_{m\geq1} {\al_{-m}\over m} z^{m} \Biggr) \eql
\sum_{m\geq0} A_{-m} z^m\,.
}
The transition matrix between the basis \eqBBac\ and \eqCBga\ is obviously
independent of 
$p$, $q$ and $q^\al$, and non-degenerate.
In the basis \eqCBga\ the matrix $\Pi^{(n)+}_{\la\mu}(p,q,q^\al)$ is
upper triangular (i.e.,
$\Pi^{(n)+}_{\la\mu}=0$ unless $\la\succeq \mu$), 
and the diagonal elements are easily computed.
We find
\eqn\eqKb{ \eqalign{
\Pi^{(n)+}(p,q,q^\al) & \eql \prod_{\la\vdash n}\ (p^{-{ {\scr1\over\scr2}}} 
q^{\al})^{\ell(\la)} \cr & \eql \prod_{{\scr r,s\geq1}\atop {\scr rs\leq n}} 
\ (p^{-{{\scr 1\over\scr 2}}} q^\al)^{p(n-rs)} \cr 
& \eql \prod_{{\scr r,s\geq1}\atop {\scr rs\leq n} }\ 
(p^{-{{\scr r\over\scr 2}}} p^{{\scr {\scr r-1\over\scr 2}}} 
q^{{{\scr s-r\over\scr 2}}} q^\al)^{p(n-rs)}\,,\cr}
}
where we have used Lemma \thCBab.  Similarly, the leading term in $q^\al$ of
$\Pi^{(n)}(p,q,pq^{-\al})$ arises from 
$\Pi^{(n)-}(p,q,q^\al)=
{\rm det}\, \Pi^{(n)-}_{\la\mu}(p,q,q^\al)$
where 
\eqn\eqCBh{
\La^-_{-\la} |\al\rangle \eql \sum_{\mu} \ 
\Pi^{(n)-}_{\la\mu}\ \al_{-\mu}^{\phantom{(}} |\al\rangle \,,
}
and is given by 
\eqn\eqKc{
\Pi^{(n)-}(p,q,pq^{-\al}) \eql 
\prod_{{\scr r,s\geq1}\atop {\scr rs\leq n}} \ 
\bigl(p^{-{ {\scr r\over\scr 2}}} p^{-{{\scr r+1\over\scr 2}}} 
q^{{{\scr r-s\over\scr 2}}} q^{\al}\bigr)^{p(n-rs)}\,.
}
This, together with the factorization \eqCBe, shows that the prefactor 
$C_n$ in \eqCBa\ is actually independent of $p$, $q$ and $q^\al$.
\Box\medskip

As a consequence of the previous proof we also have 
\thm\thCBd
\proclaim Corollary \thCBd.  
\eqn\eqCBh{
\Pi^{(n)}(p,q,q^\al) \eql D_n \ \prod_{{\scr r,s\geq1}\atop {\scr rs\leq n} }
  \left( p^{-{{\scr r\over\scr 2}}} 
  \bigl( p^{{{\scr r-1\over\scr 2}}} q^{{\scr {\scr s-r\over\scr 2}}} q^\al -
  p^{-{{\scr r-1\over\scr 2}}} q^{{{\scr r-s\over\scr 2}}} 
  q^{-\al} \bigr) \right)^{p(n-rs)} \,,
}
where $D_n$ is a constant independent of $p$, $q$ and $q^\al$.
\medskip

\proof  The Laurent polynomial 
\eqn\eqCBi{
\Pi^{\prime (n)}(p,q,q^\al) \eql \Biggl( \prod_{{\scr r,s\geq1}
\atop {\scr rs\leq n} }
  \bigl( p^{{{\scr r\over\scr 2}}} \bigr)^{p(n-rs)}\Biggr)\ 
  \Pi^{(n)}(p,q,q^\al)\,,
}
inherits the following duality invariances from $\virpq$
\eqn\eqCBj{ \eqalign{
\Pi^{\prime (n)}(p,q,q^\al) & \eql \Pi^{\prime (n)}
(p^{-1},q^{-1},q^{-\al}) \,,\cr
\Pi^{\prime (n)}(p,q,q^\al) & \eql \Pi^{\prime (n)}
(p^{-1},p^{-1}q,q^{-\al})\,.\cr}
}
This, together with the factorization \eqCBe, then uniquely determines 
$\Pi^{(n)}(p,q,q^\al)$. \Box\medskip

{}From the Kac determinant \eqCBa\ it follows that for generic $p$ and
$q$ the category $\cO$ of $\virpq$ modules is isomorphic to that 
of $\vir$.  In particular, the construction of singular vectors
[\SKAO] and a Felder type resolution of the irreducible modules 
[\LPb,\JLMP,\FJMOP] are $q$-deformations of the corresponding constructions
for $\vir$.  At $q=\sqrtN{1}$, however, the Kac determinant displays 
a large number of additional vanishing lines irrespective of the highest
weight $h$.  This indicates the existence of additional, $h$-independent,
singular vectors, and suggests that $\virpq$ at $q=\sqrtN{1}$ has a large
center.  It is important to note, though, that Corollary \thCBd\ implies
that the homomorphism $\imath:
M(h(\al)) \rightarrow F(\al)$, for generic $p$ and $\al$, is a bijection
even for $q=\sqrtN{1}$. 
In the following section we use this fact to establish results for
$\virpq$ at $q=\sqrtN{1}$ by using the free field realization.

%
\newsec{The center and representations of $\virpq$ at roots of unity}

In this section we analyze $\virpq$ and its representations for
$q$ a primitive $N$-th root of unity.  To illustrate the main 
features we first discuss the simplest case, $N=2$, before we proceed
to general $N$.  In the last part of this section we make the analysis
more explicit in the limit $t\to\infty$.

\subsec{Representations of $\virpq$ for $q=-1$}

For $q=-1$ we have 
\eqn\eqDAa{
f(x) \eql { 1+x\over 1-x}\,,
}
i.e., $f_0=1$ and $f_l=2$ for all $l\geq1$.  Thus $\virpq$ at $q=-1$ is
given by
\eqn\eqDAb{
[T_m,T_n] + 2 \sum_{l>0} \left(T_{m-l}T_{n+l} - T_{n-l}T_{m+l} \right)
  \eql c_m \de_{m+n,0}\,,
}
where (cf.\ \eqBAc)
\eqn\eqDAc{
c_m \eql -2  \left( {1+p\over 1-p}\right) (p^m - p^{-m}) \,.
}
Note that, for all $m$ and $n$, the sum over $l$ in \eqDAb\ is actually 
finite.
\thm\thDAa
\proclaim Theorem \thDAa.  For $q=-1$, $\virpq$ is equivalent to the algebra
defined by the relations
\eqn\eqDAd{
\{T_m,T_n\} \eql \cases{
  2 (-1)^{{\scr m-n\over\scr 2}} T_{{\scr m+n\over\scr 2}}
T_{{\scr m+n\over\scr2}} + 
  \widetilde{c}_m \de_{m+n,0} & {for $m+n\in2\ZZ\,,$} \cr
  0 & {otherwise}$\,,$ \cr}
}
where
\eqn\eqDAe{
\widetilde{c}_m \eql 2 \left( p^{{\scr m\over\scr2}} - p^{-{\scr m\over\scr2}} 
\right)^2\,.
}
\medskip

\proof The proof is straightforward if one uses the identity
\eqn\eqDAf{
\sum_{l=0}^m (-1)^l f_l c_{m-l} \eql \widetilde{c}_m\,,
}
which is proved by induction. \Box\medskip

\ni {\bf Remark.} It is possible to arrive at the commutators \eqDAd\
directly, using the free field realization, by exploiting a different 
factorization of the exchange matrix $S_{TT}(x) = f(x)^{-1}f(1/x) =
-1$, namely $S_{TT}(x) = g_+(x)^{-1} g_-(1/x)$ with $g_+(x) = -
g_-(x) =1$.  See [\FRb] for details regarding this procedure.\medskip

Note that the commutation relations \eqDAb\ can be 
considered as an equation for the symmetrization of the product of two
generators.  In Section 4.3 we generalize this, in the limit $t\to\infty$, 
to the symmetrization of a product of $N$ generators for arbitrary $q\in\CC$.
\medskip

The next result is 
a simple consequence of the commutation relations \eqDAb.
\thm\thDAb
\proclaim Theorem \thDAb. The elements $(T_n)^2$, $n\in\ZZ$, are in the center
of $\virpq$ at $q=-1$.  In particular, this implies that
the vectors $(T_{-n})^2|\hw\rangle$ are
singular for all $p, \hw\in \CC$ and $n\in\NN$.
\medskip

\ni {\bf Remark.} On the basis \eqCAf\ of $M(\hw)$, the 
action of $T_0$ is upper-triangular, i.e.,
\eqn\eqDAg{
T_0 |\la;\hw\rangle \eql (-1)^{\ell(\la)} \hw |\la;\hw\rangle + 
 \sum_{ \mu < \la} c_{\la\mu} | \mu;\hw\rangle\,,
}
e.g.,
at level $2$ on the basis $\{T_{-1}T_{-1}|\hw\rangle,T_{-2}|\hw\rangle \}$
we have 
\eqn\eqDAh{
T_0 \eql \left( \matrix{ \hw & -2 \cr 0  & -\hw \cr} \right)\,.
}
In particular note that
$T_0$ in \eqDAh\ is not diagonalizable for $\hw=0$.

Let $M'(\hw)$ be the module obtained from
$M(\hw)$ by dividing out the submodule generated by the singular vectors 
$(T_{-n})^2|\hw\rangle, n\geq1$.
\thm\thDAca
\proclaim Theorem \thDAca.  The Kac determinant $G^{\prime (n)}$ of 
$M'(\hw)_{(n)}$ is given by 
\eqn\eqDAla{
G^{\prime (n)} \eql C_n  \prod_{{\scr r\geq1}\atop {\scr r\leq n} }
  (\hw^2 - \hw^2_{r,1})^{q_2(r;n-r)} \,,
}
where $C_n$ is a constant independent of $p$ and $h$, and the integers 
$q_2(r;n)$ are determined by
\eqn\eqDAlb{
\prod_{{\scr n\geq1}\atop {\scr n\neq r}} 
\ (1+x^n) \eql \sum_{n\geq0} \ q_2(r;n) \,x^n\,.
}\medskip

\proof Observe that 
\eqn\eqDAk{
\widetilde{c}_m \eql 2 \left( p^{{\scr m\over\scr 2}} - 
p^{-{\scr m\over\scr 2}} 
\right)^2 \eql - 2 (-1)^m \hw^2_{m,1} \,,
}
while the character of the submodule of $M'(\hw)$ generated by 
$T_{-r}|\hw\rangle$ is given by \eqDAlb. \Box\medskip

\thm\thDAc
\proclaim Theorem \thDAc. Let $M'(\hw)$ be the quotient module
defined as above.
\item{(i)} $M'(\hw)$ is irreducible provided $\hw^2 \neq 
  \hw_{r,1}^2$ for all $r\geq1$.
\item{(ii)} The character of $M'(\hw)$ is given by
\eqn\eqDAi{
{\rm ch}_{M'} (x) \eql \prod_{n\geq1} (1 + x^n) \eql 
  \prod_{n\geq1}\left(\frac{1-x^{2n}}{1-x^{n}}\right) \,.
}
\item{(iii)} The `Witten index' of $M'(\hw)$ is given by
\eqn\eqDAj{
{\rm Tr} (T_0 x^d) 
  \eql \hw \prod_{n\geq1} (1 - x^n)\,. 
}
\medskip 

\proof 
\item{(i)} Follows from Theorem \thDAca.
\item{(ii)} In $M'(\hw)$ we have an (orthogonal) basis of monomials 
\eqn\eqDAm{
T_{-\la_1} \ldots T_{-\la_\ell} |\hw\rangle, \quad \la_1>\ldots >\la_\ell>0\,,
}
indexed by partitions with no equal parts.
\item{(iii)} Observe that, in $M'(\hw)$,
\eqn\eqDAn{
T_0 \left( T_{-\la_1} \ldots T_{-\la_\ell} |\hw\rangle \right) \eql
  (-1)^{\ell(\la)} \hw \left( T_{-\la_1} \ldots T_{-\la_\ell}
 |\hw\rangle \right)\,.
}
This concludes the proof. \Box\medskip

\ni {\bf Remark.}  Clearly, the module $M'(\hw)$ can be realized in
terms of free (Ramond) fermions, and is isomorphic to the irreducible 
$\vir$ module at $c={1\over2}$ and $\De={1\over16}$.  We thus have
an  action of $\virpq$, for $q=-1$, on a Virasoro minimal model module.
\medskip

Now, if $\hw=\pm \hw_{m,1}$ for some $m$, then $M'(\hw)$ is
reducible.  We leave the general analysis of this situation for further
study.  Here, let us just remark, that if in addition $t$ is a
primitive $M$-th root of unity%
\foot{Note that in this case, $p$ is a root of $-1$ and thus   
the algebra $\virpq$ is not uniquely
defined (cf.\ the remark in Section 2.1).  
The particular algebra we are discussing here
corresponds to first taking $q\to-1$ and then $t\to \root M\of 1$, i.e.,
to the solution \eqDAa\ of the recursion relations \eqBAg.},
e.g., $t=\exp( {2\pi i\over M})$, then if the equation 
\eqn\eqDAo{
\hw \eql \pm \hw_{m,1} = \pm 2 \cos( \pi ({m\over M} -{1\over2})) 
  \eql \pm 2 \sin(  {\pi m\over M})\,,
}
holds for one particular $m=m_0$, it holds for all 
$m= m_0\, ({\rm mod}\, M)$
and $m= M-m_0\, ({\rm mod}\, M)$.
It is straightforward to verify the following theorem.
\thm\thDAd
\proclaim Theorem \thDAd. Suppose $t=\exp( {2\pi i\over M})$ and 
$\hw = \pm 2 \sin(  {\pi m_0\over M})$ for some $m_0\in\{0,1,\ldots,M-1\}$.
Then the vectors $T_{-m}|\hw\rangle$ are singular in $M'(\hw)$ for
all $m= m_0\,({\rm mod}\, M)$ and $m= M-m_0\,({\rm mod}\,M)$.  
Let $M''(\hw)$ be the module obtained
from $M'(\hw)$ by dividing out the ideal generated by these singular
vectors.  The module $M''(\hw)$ is irreducible and has character
\eqn\eqDApaa{
{\rm ch}_{M''} (x) \eql \prod_{n\geq1} 
  \left({{1+x^n\over 1+x^{m_0+nM}}}\right)\,,
}
for $m_0 =0$ or $m_0={M/2}$, and 
\eqn\eqDApab{
{\rm ch}_{M''} (x) \eql  \prod_{n\geq1} 
  \left({{1+x^n\over (1+x^{m_0+nM})(1+x^{(M-m_0)+nM})}}\right)\,,
}
otherwise.
\medskip

\ni {\bf Remark.} The irreducibility of $M''(h)$ for, e.g., $M=3$ and
$h=h_{1,1}$ is in apparent conflict with case (ii) of Theorem \thCc.
However, this is precisely a manifestation of the fact that the 
algebra is not uniquely defined for $p$ a root of $-1$.\medskip

\ni {\bf Remark.} Note that the vectors $T_{-m}|\hw\rangle$ do not have 
to be singular in $M(\hw)$.  For example, consider $M=3$, $\hw=0$ ($m_0=0$).
Then 
\eqn\eqDAq{
T_1  T_{-3}|\hw\rangle \eql 2 T_{-1} T_{-1} |\hw\rangle\,,
}
which means that 
$T_{-3}|\hw\rangle$ is primitive, but non-singular, in $M(\hw)$.
As a consequence, the Verma module $M(\hw=0)$ for $q^2=t^3=1$ has a submodule
that is not generated by singular vectors, namely, the submodule generated
by $T_{-3}|\hw\rangle$ and $T_{-1} T_{-1} |\hw\rangle$.  This situation 
does not occur for Verma modules of the Virasoro algebra, but is common
in the case of $\cW$-algebras [\BMP].\medskip

\subsec{Representations of $\virpq$ for $q=\sqrtN{1}$}

In this section we consider the case where $q$ is an arbitrary
primitive $N$-th root of unity with $N>2$. Our main goal is to construct
explicitly the series of singular vectors in the Verma module
$M(h)$ for generic $p$, 
that are indpendent of $h$, and to analyze the structure of
the resulting quotient module $M'(h)$. The results, given in Theorems
\thDBb\ and \thDCfa, generalize those from the previous section.
However, unlike for $q=-1$, the commutation relations of $\virpq$ are
extremely cumbersome for
$q=\sqrtN{1},\ N>2$, and it is difficult to
establish any of these results directly. For that reason we 
resort to the free field realization of Section 2.3, which
considerably simplifies the entire analysis, as we now show.
\thm\thDBxa
\proclaim Lemma \thDBxa. For $q=\sqrtN{1}$, the oscillators $\al_m$ with
$m= 0\, {\rm mod}\,N$ generate the center of the $q$-Heisenberg
algebra, $\cH_{p,q}$. In particular, they therefore also commute with
$\imath(T(z))$. \medskip

\proof  This result is an obvious consequence of the commutation relations
\eqBBa\ in $\cH_{p,q}$ for $q=\sqrtN{1}$. \Box\medskip

Since polynomials in the oscillators $\al_{-nN}$, $n\in\NN$,
commute with $\imath(T(z))$, upon acting on the vacuum
they give rise to singular vectors in the
Fock space $F(\al)$.  Given that
$\imath: M(h)\rightarrow F(\al)$ is an isomorphism (cf.\ Section 3), we
conclude that there are corresponding singular vectors in the Verma
module $M(h)$, that are independent of a particular value of $h$. We
will now proceed to compute those vectors explicitly.
\thm\thDBa
\proclaim Lemma \thDBa.  For $q=\sqrtN{1}$ and $p$ generic 
\eqn\eqDBe{ \eqalign{
\lim_{z_i\rightarrow q^{N-i}z}\ \Biggl(\,
 \prod_{i<j}f({z_j\over z_i})\Biggr)& \,
\imath( T(z_1)\ldots T(z_N) ) \cr & \eql  \NO{\La^+(zq^{N-1})\ldots \La^+(z)} 
 + \NO{\La^-(zq^{N-1})\ldots \La^-(z)}\,.\cr}
}
\medskip

\proof  Following the steps in the proof of Theorem \thBBa, we find
\eqn\eqDBf{
\Biggl(\, \prod_{i<j}f({z_j\over z_i})\Biggr)\,
\imath( T(z_1)\ldots T(z_N)) \eql \sum_{\ep_1,\ldots,\ep_N}\ 
\Biggl(\,\prod_{i<j}
F^{\ep_i\ep_j}({z_j\over z_i})\Biggr) \  
\NO{\La^{\ep_1}(z_1)\ldots \La^{\ep_N}(z_N)}\,.
}
For generic $p$, $F^{\ep_i\ep_j}(x)$ has a
well-defined limit for $x\rightarrow q^k$ as long as $q^k\neq1$.  It
follows that the right hand side of \eqDBf\ has a
well-defined limit for $z_i\rightarrow zq^{N-i}$. In this limit we find
(cf.\ [\FRb], Section 6.1)
\eqn\eqDBg{ \eqalign{
{\rm \eqDBf} \quad\rightarrow\quad &
\NO{\La^+(zq^{N-1})\ldots \La^+(z) } +
\NO{\La^-(zq^{N-1})\ldots \La^-(z) } \cr
& + \sum_{k=1}^{N-1} b_{k,N}(q) \ 
\NO{\La^+(zq^{N-1})\ldots \La^+(zq^{N-k})\La^-(zq^{N-(k+1)})\ldots
\La^-(z)}\,,\cr}
}
where
\eqn\eqDBh{ \eqalign{
b_{k,N}(q) & \eql \prod_{i=1}^k\,\prod_{j=k+1}^N\ F^{+-}(q^{i-j}) \cr
& \eql \prod_{i=1}^k\,\prod_{j=k+1}^N\ 
{(1-q^{i-j-1})(1-p^{-1}q^{i-j+1})
\over (1-q^{i-j})(1-p^{-1}q^{i-j})}\,.\cr}
}
To establish \eqDBg, we first consider terms in \eqDBf\ that have
$\La^-(z_i)$ to the left of some $\La^+(z_j)$. Choosing the rightmost
such $\La^-(z_i)$ we see that there is a factor
$$
\ldots F^{-+}({z_{i+1}\over z_i})\ldots \La^-(z_i)\La^+(z_{i+1}) \ldots\,,
$$
which vanishes in the limit $z_i\to zq^{N-i}$ because $F^{-+}(q^{-1})=0$. The
remaining terms yield \eqDBg.
Finally, we note that for $q=\sqrtN{1}$, we have 
\eqn\eqKd{
b_{k,N}(q) \eql 0\,,\qquad k=1,\ldots,N-1\,,
}
because of the factor $(1-q^{-N})$ that arises by setting
$i=1$ and $j=N$ in \eqDBh. \Box\medskip

Now, it follows from the explicit expression \eqBBd\ 
that for $q=\sqrtN{1}$, the terms 
$$
\NO{\La^\pm(zq^{N-1})\ldots \La^\pm(z)}
$$
in \eqDBe\ have an expansion in terms of $\{\al_{nN}, n\in\ZZ\}$. 
This observation, together with Lemma \thDBa, yields the following theorem.
\thm\thDBb
\proclaim Theorem \thDBb.  For $q=\sqrtN{1}$ and $p$ generic
\eqn\eqDBi{
\Psi(z) \eql \lim_{z_i\rightarrow zq^{N-i}}\ 
\Biggl( \prod_{i<j}f({z_j\over z_i})\Biggr) \,
T(z_1)\ldots T(z_N) |\hw\rangle\,,
}
is a well-defined generating series of singular vectors in $M(\hw)$.
\medskip

\proof The result of Lemma \thDBa\ is that $\imath(\Psi(z))$ is a
well-defined series of singular vectors in $F(\al)$. Since $\imath$ is
an isomorphism, this also proves the theorem. \Box\medskip

Although \eqDBi\ defines $\Psi(z)$ as a series in the products of
modes of $T(z)$ acting on the vacuum, interpreting this result
directly within the Verma module must be done with some caution. On the
one hand, upon expanding $T(z_i)$ into power series, we find that the
resulting modes $T_m$ have both $m\leq 0$ and $m>0$. In the latter
case one can commute those modes to the right until they annihilate
the vacuum. However, as one can see from \eqBAa, this yields
additional infinite summations.  Therefore, the product $T(zq^{N-1})\ldots
T(z)$ turns out to be a divergent series, even when acting on the
vacuum.

On the other hand, it follows directly from the second expression in
\eqBAb\ that, for $q=\sqrtN{1}$, 
\eqn\eqDBr{ 
f(x)f(xq)\ldots f(xq^{N-1}) \eql 1\,.  
} 
Since the product of $f(z_j/z_i)$ in \eqDBi\
has a factor of $f(q)\ldots f(q^{N-1})= f(1)^{-1}=0$ in the limit, it 
vanishes.\medskip

\ni {\bf Remark:} Note that, by using \eqDBr\ and the exchange relations
\eqBBexch, the operator $T(zq^{N-1})\ldots T(zq)T(z)$ is formally in
the center of $\virpq$.  However, as we have seen above, it diverges.
The products of $f(z_j/z_i)$ in \eqDBi\ provide a convenient
regularizing factor (see, also [\FRb]).\medskip

Now, let us carefully expand \eqDBi\ in modes. We have
\eqn\eqDBia{\eqalign{
\Biggl( \prod_{i<j}  & f({z_j\over z_i})\Biggr) \,
T(z_1)\ldots T(z_N) |\hw\rangle 
\eql \sum_{m_1,\ldots,m_N\in\ZZ}\, z_1^{m_1}
\ldots z_N^{m_N}\,
\sum_{l_{ij}\geq0} \ \Biggl( \prod_{i<j} f_{l_{ij}} \Biggr) \cr
& \times T_{-m_1-l_{12}-\ldots-l_{1N}} 
 T_{-m_2+l_{12}-l_{23}-\ldots -l_{2N}}
\ldots T_{-m_N+l_{1N}+\ldots+l_{N-1 N}} |\hw\rangle\,. \cr}
}
Upon introducing 
Littlewood's raising operators $R_{ij},\ i<j$ [\Lit,\MD], acting 
on monomials as
\eqn\eqDBl{
R_{ij}\, T_{-m_1}\ldots T_{-m_N} |\hw\rangle \eql
T_{-m_1}\ldots T_{-(m_i+1)}\ldots T_{-(m_j-1)}\ldots  
T_{-m_N} |\hw\rangle \,,
}
we can write \eqDBia\ more succinctly,
\eqn\eqDBib{
{\rm \eqDBia} ~\rightarrow~ \sum_{m_1,\ldots,m_N\in\ZZ}\, z_1^{m_1}
\ldots z_N^{m_N}\, 
\Biggl(\,\prod_{i<j} f(R_{ij})\Biggr) \ 
T_{-m_1} \ldots T_{-m_N} |\hw\rangle\,,
}
where 
\eqn\eqDBk{
f(R_{ij}) \eql \sum_{l\geq0} \ f_l\, (R_{ij})^l\,.
}
Note that in this notation the commutation relations \eqBAa\ simply
become
\eqn\eqDBcrx{ 
f(R_{12})T_{m}T_{n}\eql f(R_{12})T_{n}T_{m} +c_m \,\de_{m+n,0}\,.
}

Now consider \eqDBib. Clearly, at a given level in the Verma module,
say level $d$, only the terms satisfying $\sum m_i=d$ contribute.
There are only a finite number of such terms that have all $m_i\geq0$. In
all other terms we can commute the $T_{-m_i}$ with $m_i<0$ to the
right, where they annihilate the vacuum. In fact, this would have been
quite simple if there was no central term in \eqDBcrx, as in this case
there would be a full symmetry in $m_1,\ldots,m_N$. However, because
of the central term in the commutation relations, we do obtain
subleading terms, i.e., with the product of a smaller number of the
$T_{m_i}$.  Moreover, those central terms result in infinite sums over
$m_i<0$. It is easy to see that for $|p|$ sufficiently small, one can
sum up those series, and the final expression is manifestly
well-defined in the limit $z_i\rightarrow q^{N-i}z$.

To illustrate this procedure let us consider the simpler case with
$N=2$. Here we find
\eqn\eqDBntx{
\sum_{m_1,m_2\in\ZZ} z_1^{m_1}z_2^{m_2} f(R_{12}) 
T_{-m_1}T_{-m_2}|h\rangle \eql
\left(\sum_{\la}m_{\la}(z_1,z_2)f(R_{12}) T_{-\la_1}T_{-\la_2}
+ c({z_2\over z_1})\right) |h\rangle\,,
}
where $m_\la(z_1,z_2)$ are the monomial symmetric polynomials and the sum is
over all partitions.

{}For $N>2$ it becomes considerably more difficult to carry out
this calculation. In Appendix B we summarize the result for $N=3$. 
Here let us only say that, as follows from the discussion above, the 
leading term in \eqDBi, i.e., the term with
$N$ generators $T_n$, is given by
\eqn\eqDBj{
\sum_{{\scr \la}\atop {\scr \ell(\la)\leq N}} \ 
m_\la(z,zq,\ldots,zq^{N-1})\ \Biggr(\prod_{i<j} f(R_{ij})
\Biggl) \ T_{-\la_1}\ldots T_{-\la_N} |\hw\rangle\,,
}
where the sum is over all partitions $\la$.

In particular, it follows from \eqDBj\ that the singular vectors $\Psi_{-d}
|\hw\rangle$ for $d=mN,\ m\in \NN$, are non-vanishing and independent, i.e., 
none of the $\Psi_{-d}|\hw\rangle$ is in the submodule generated by the
other ones.  The vanishing of the leading term for $d\neq
0{\,\rm mod}\,N$ follows from part (i) of the following lemma.
\thm\thDBd
\proclaim Lemma \thDBd.  For $q=\sqrtN{1}$ we have 
\item{(i)}
\eqn\eqKe{ \eqalign{
m_\la & (z,zq,\ldots,zq^{N-1}) \cr & \eql \cases{
q^{{1\over2}mN(N-1)} \ 
(K(1)^{-1}K(q))_{\la,(m^N)}\ z^{|\la|} & if $\la\vdash mN$ for
some $m\in\ZZ_{\geq0}\,,$\cr
0 & otherwise$\,.$\cr} \cr}
}
\item{(ii)}
\eqn\eqKf{
m_{(\la_1+n,\ldots,\la_N+n)}(z,zq,\ldots,zq^{N-1}) \eql 
m_{(\la_1,\ldots,\la_N)}(z,zq,\ldots,zq^{N-1})\,,\qquad \forall n\in\NN\,.
}
\medskip

\proof For (i), use the following expansion of $m_\la(x)$ in terms of 
Hall-Littlewood polynomials
which immediately follows from \appAxxv\ and \appAzx\ in Appendix A,
\eqn\eqDBo{
m_\la(x) \eql \sum_\mu \ (K(1)^{-1}K(q))_{\la\mu}\ P_\mu(x;q)\,.
}
Then use (cf.\ [\MD])
\eqn\eqDBp{
P_\la(z,zq,\ldots,zq^{N-1};q)=q^{n(\la)}\ \bin{N}{m(\la)}\ z^{|\la|}\,,
}
where 
\eqn\eqDBq{
n(\la)\eql \sum_{i\geq 1}(i-1)\la_i\,,
}
and 
\eqn\eqLb{
\bin{N}{m(\la)} ~\equiv~ \bin{N}{m_1(\la) \ldots m_N(\la)} \eql
{ (q)_N \over (q)_{m_1(\la)} \ldots (q)_{m_N(\la)} }\,.
}
Now, for $q=\sqrtN{1}$, the $q$-multinomial \eqLb\ vanishes unless 
$m_i(\la)=N$ for some $i\in\ZZ_{\geq0}$, i.e., only the terms for which
$\mu=(m^N)$ for some $m\in\ZZ_{\geq0}$ contribute in \eqDBo.
This proves (i).  Part (ii) is obvious from the definition of $m_\la(x)$.
\Box\medskip

We have not succeeded in finding a more explicit, but still tractable,
expression for the singular vectors of \eqDBi\ for arbitrary $N$.  From 
the explicit expressions of some singular vectors at $q=\sqrtN{1},
\,N=3,4$, in Appendix B, it is however clear that the expressions drastically
simplify in the limit $t\to\infty$ (or, equivalently, $t\to0$).  Indeed,
in all examples only the leading term $(T_{-n})^N|h\rangle$ survives.
In section 4.3 we analyze this limit in more detail.\medskip

For $q=\sqrtN{1}$,
let $\cH'_{p,q}$ denote the reduced $q$-Heisenberg algebra, i.e.,
$\cH'_{p,q}$ with the oscillators $\{ \al_n\, |\, n= 0\, {\rm mod}\,N\}$
removed, and denote by $F'(\al)$ the Fock space of $\cH'_{p,q}$.
\thm\thDBe
\proclaim Theorem \thDBe. For $q=\sqrtN{1}$, and 
generic $p$, we have a realization 
of $\virpq$ on the sub Fock space $F'(\al) \subset F(\al)$.  
This realization is irreducible for generic $\al,p \in \CC$
and the character is given by 
\eqn\eqDBa{
{\rm ch}_{F'} (x) \eql \prod_{n\geq1}\left({1-x^{nN} \over
1-x^{n}}\right) 
\eql \sum_{n\geq 0} \ p_N(n) \,x^n\,,
}
where $p_N(n)$ is the number of partitions of $n$ with parts 
not equal to a multiple of $N$.
\medskip

The proof of this theorem is clear except for the irreducibility of $F'(\al)$.
This will follow from the result of Theorem \thDCfa.

Let $M'(\hw)$ denote the module obtained from $M(\hw)$ by dividing 
out the submodule generated by the singular vectors of \eqDBi.  
To investigate the irreducibility of $M'(\hw)$
we make use of the following
\thm\thDCf
\proclaim Theorem \thDCf. The Kac determinant of the module $M'(\hw)$
is given by 
\eqn\eqDCo{
\widetilde{G}^{\prime (n)} \eql C_n \prod_{ {\scr r\geq1, 
 1\leq s\leq N-1} \atop {\scr rs\leq n}}
 (\hw^2 - \hw_{r,s}^2)^{q_N(r,s;n-rs)}
 \prod_{{\scr r,s\geq 1,\ rs\leq n}\atop {\scr r\neq 0\,{\rm mod}\,N}}
  \left( {1-t^{-r}\over 1+p^r} \right)^{p_N(n-rs)} \,,
}
where $C_n$ is a constant independent of $p$ and $h$, $p_N(n)$
is defined in \eqDBa, and 
the integers $q_N(r,s;n)$ are determined by 
\eqn\eqDCoa{
(1+x^r+\ldots+x^{r(N-1-s)})\,
\prod_{{\scr n\geq1 }\atop {\scr n\neq r}}
 (1 + x^n + \ldots + x^{n(N-1)}) \eql
\sum_{n\geq0} \ q_N(r,s;n)\,x^n\,.
}
\medskip

\proof The proof is completely analogous to the proof
of Theorem \thCBa.  In particular, the second term arises
from the Kac determinant $g^{\prime(n)}$ of $F'(\al)$ while
the first term is a remnant of the vanishing lines \eqCBf.
Note that $h_{r,s} = \pm h_{r,s+N}$ for $q=\sqrtN{1}$. \Box\medskip

The generalization of Theorem \thDAc\ reads
\thm\thDCfa
\proclaim Theorem \thDCfa. Let $q=\sqrtN{1}$ and 
let $M'(\hw)$ be defined as above, then
\item{(i)} $M'(\hw)$ is irreducible provided $h^2\neq h^2_{r,s}$ for all
$r\geq1$ and $1\leq s\leq N-1$.
\item{(ii)} The character of $M'(\hw)$ is given by
\eqn\eqDCn{
{{\rm ch}}_{M'}(x) \eql \prod_{n\geq1} (1+ x^n + \ldots + x^{n(N-1)})\,.
}
\medskip

\ni {\it Proof of Theorem \thDBe:}  It remains to prove the irreducibility of
$F'(\al)$ for generic $\al$ and $p$.  This follows from the irreducibility
of $M'(\hw(\al))$ (Theorem \thDCfa\ (i)) and 
the equality of characters
\eqn\eqDCnz{
\prod_{n\geq1}\left({1-x^{nN} \over 1-x^{n}}\right) 
\eql \prod_{n\geq1} (1 + x^n + \ldots + x^{n(N-1)})\,,
}
which imply that $\imath: M'(\hw(\al)) \rightarrow F'(\al)$ is an isomorphism.
\Box\medskip

It is useful to have an explicit expression for the image
of a Verma module monomial $T_{-\la}|\hw\rangle$ under the map
$\imath: M(h(\al))\rightarrow F(\al)$.  To describe the result, let us 
identify $F(\al)$ with the ring of symmetric functions 
$\La \otimes \QQ(p,q)$ over $\QQ(p,q)$
through the isomorphism $\jmath$,
\eqn\eqDBsa{
\jmath (\al_{-\la_1}\ldots \al_{-\la_n} |\al\rangle) \eql p_\la\,,
}
where $p_\la$  are the power sum symmetric 
functions (see, Appendix A).  We then have 
\thm\thDBc
\proclaim Theorem \thDBc.  Let $\la$ be a partition of length $\ell(\la)=n$,
then 
\eqn\eqDBma{
\imath(T_{-\la_1}\ldots T_{-\la_n} |\hw\rangle) \eql
\sum_{\ep_1,\ldots,\ep_n} \ 
\La^{\ep_1}_{-\la_1} \ldots \La^{\ep_n}_{-\la_n} |\al\rangle \,,
}
while
\eqn\eqDBm{
\La^{\ep_1}_{-\la_1} \ldots \La^{\ep_n}_{-\la_n} |\al\rangle \eql
(p^{-{1\over2}} q^\al)^{\sum \ep_i} \ 
\Biggl( \prod_{i<j} \ f(R_{ij})^{-\ep_i\ep_j} \Biggr)
\ h^{\ep_1}_{\la_1}(x)\ldots
h^{\ep_n}_{\la_n}(x) \,,
}
where
\eqn\eqDBn{
h_{n}^+(x) \eql h_n(x)\,,\qquad 
h_{n}^-(x) \eql p^{-n} e_n(-x)\,.
}
and $f(x)$ as in \eqBAb.
\medskip

We recall that $h_n(x)$ and $e_n(x)$ are, respectively, the completely
symmetric functions and the elementary symmetric functions.  
[The expression \eqDBm\ is to be understood as follows:
first replace $h_{n_i}^{\ep_i}(x)$ by the appropriate expression in
\eqDBn, then act with the $f(R_{ij})$ on any combination of $h_n$'s
and $e_n$'s as in \eqDBl.  In particular, the $f(R_{ij})$ act only on
the subscript of the symmetric function, and {\it not} on the
prefactor $p^{-n}$.]\medskip

\proof As in [\Jib], Proposition 3.9. \Box\medskip

As an application  of Theorem \thDBc, consider the leading 
order term \eqDBj\ of the singular vector $\Psi(z)$.  Applying the 
map $\imath$ we see that the term of leading order in $q^\al$ (i.e.,
the term of order $\cO( (q^\al)^N )$) is proportional to
\eqn\eqKg{ \eqalign{
\sum_{{\scr \la}\atop {\scr \ell(\la)=N}}
 \ m_{\la}(z,zq,\ldots,zq^{N-1}) h_\la(x)
& \eql \sum_{m_1,\ldots,m_N\in\ZZ_{\geq0}}\ 
h_{m_1}(x) h_{m_2}(xq) \ldots h_{m_N}(xq^{N-1})\cr & \eql
\sum_{m\geq0}\ h_m(x^N) z^{mN}\,,\cr}
}
where, in the last step, we have used 
\eqn\eqKh{
\prod_{i=1}^{N} \ H(xq^{i-1};t) \eql H(x^N;t^N)\,.
}
Now,
\eqn\eqKi{
h_m(x^N) \eql \sum_{\la\vdash m}\ z_\la^{-1}\, p_\la(x^N) \eql 
 \sum_{\la\vdash m}\ z_\la^{-1}\, p_{N\la}(x)\,,
}
so that indeed (cf.\ \eqDBsa) the vector \eqDBj\ gives rise to
a singular vector at the leading order in $q^\al$.  A similar 
computation also shows that the leading order in $q^{-\al}$ 
works out.  The subleading terms in \eqDBi\ are required to make
the subleading orders in $q^\al$ work.

As a second application of Theorem \thDBc\ one can verify that,
for $q=-1$, the vectors $(T_{-n})^2|\hw\rangle$ are indeed
singular (cf.\ Theorem \thDAb).\medskip

\subsec{Representations of $\virpq$ for $q=\sqrtN{1}$ in 
 the limit $t\to\infty$}

In the previous section we have analyzed $\virpq$ and its
representations for $q=\sqrtN{1}$ and generic $p$.  In this section we
will make the analysis even more explicit in the limit $t\to\infty$
(or, equivalently, the limit $t\to0$), where, as we have emphasized
already, we expect some drastic simplifications.
At the same time the limit $t\to\infty$ is a generic
point, in the sense that the structure of the algebra and the modules
is similar to that at other generic values of $t$, e.g.,
multiplicities of singular vectors at $t=\infty$ are equal to the
multiplicities at a generic $t$.  It turns out that the most efficient
method of studying this case is to exploit an interesting connection
to Hall-Littlewood polynomials, whose basic properties have been
summarized in Appendix A.
\medskip

The deformed Virasoro algebra $\virpq$ at $t=\infty$,  denoted by
$\tvirq$, is generated by
$\{\wT_n,\,n\in\ZZ\}$, where 
\eqn\eqDCb{
\wT_m \eql \lim_{t\rightarrow\infty} T_m \, p^{{|m|\over2}}\,,
} 
satisfy the relations given by the following theorem (cf.\ [\AKOSb]):
\thm\thDCa
\proclaim Theorem \thDCa. In the limit $t\to\infty$ we have 
\eqn\eqDCaa{ 
[\widetilde{T}_m,\widetilde{T}_n]_q  \eql (q-q^{-1}) \sum_{l\geq1}\, q^{l}\, 
\widetilde{T}_{n-l}\widetilde{T}_{m+l} + (1-q)  \de_{m+n,0}, 
\qquad {\rm if}\quad m>0>n\,,
}
\eqn\eqDCab{ 
[\widetilde{T}_m,\widetilde{T}_n]_q  \eql -(1-q) \sum_{l=1}^{m-n-1}\,
\widetilde{T}_{m-l}\widetilde{T}_{n+l},
\qquad {\rm if}\quad m>n>0 \quad {\rm or} \quad 0> m>n\,,
}
\eqn\eqDCac{
[\widetilde{T}_0,\widetilde{T}_m]_q  \eql -(1-q) \sum_{l=1}^{-m-1} \,
\widetilde{T}_{-l}
\widetilde{T}_{m+l} + (q-q^{-1}) \sum_{l\geq1} q^{l}\, \widetilde{T}_{m-l}
\widetilde{T}_l, \qquad {\rm if}\quad 0>m\,,
}
\eqn\eqDCad{
[\widetilde{T}_m,\widetilde{T}_0]_q  \eql -(1-q) \sum_{l=1}^{m-1} 
\widetilde{T}_{m-l}
\widetilde{T}_l + (q-q^{-1}) \sum_{l\geq1} q^{l} \widetilde{T}_{-l}
\widetilde{T}_{m+l}, \qquad {\rm if}\quad m>0\,,
}
where $ [-,-]_q$ is the $q$-commutator, 
\eqn\eqDCc{
[x,y]_q \eql xy - qyx\,.
}
\medskip

\thm\thDCb

\proof By a direct expansion of the commutation relation \eqBAa\ to
the leading order in $t$ using Lemma 
\thDCb\ below. \Box
\medskip

\proclaim Lemma \thDCb. For $t\to\infty$ we have 
\eqn\eqDCd{ \eqalign{
f_0 & \eql 1,\cr
f_l & \eql (1-q) + \cO(t^{-1})\,,\cr}
}
and
\eqn\eqDCe{
f_l - f_{l+m} \eql - q^{2l-1} (1-q^2) t^{-l} + \cO(t^{-l-1})\,,
}
for all $m>l\geq1$.
\medskip

\ni
{\bf Remark.} Note that \eqDCd\ is equivalent to 
\eqn\eqDCda{
f(x) \eql {1-qx\over 1-x} + \cO(t^{-1})\,.
}
\medskip

\proof By an explicit expansion of \eqBAb. \Box
\medskip

Let us point out some obvious simplifications of the relations
\eqDCaa\ -- \eqDCad\ over \eqBAa. First, unlike in $\virpq$, we have 
subalgebras $\tvirq^\pm$ of $\tvirq$ generated by the $\wT_m$ with
$\pm m>0$. The relation \eqDCab\ that defines those subalgebras has
only a finite number of terms on the right hand side. Moreover,
\eqDCab\ is invariant under the `shift transformation'
induced by $(m,n)\rightarrow (m+k,n+k)$, $k\in\ZZ$, as long as one
remains within a given subalgebra.  Secondly, the commutation
relations between the positive and negative mode generators do not
produce the $\widetilde T_0$, and the right hand side of \eqDCaa\ is
already in  ordered form. Finally, the form of \eqDCac\ and
\eqDCad\ suggests that one should be able to represent $\widetilde
T_0$ as a sum $\widetilde T_0=\widetilde T_0^++\widetilde T_0^-$,
where $\widetilde T_0^\pm$ extend $\tvirq^\pm$ by a zero mode
generator.
\medskip

Now, let us consider the Verma module
$\wM(\hw) = \amalg_{n\geq0} \wM(\hw)_{(n)}$, in which we introduce a
basis
\eqn\eqDCj{
|\la;\hw\rangle \eql \wT_{-\la_1} \ldots \wT_{-\la_\ell} |\hw\rangle\,,
 \qquad \la\vdash n\,, } of $\wM(\hw)_{(n)}$ as in \eqCAf.
\thm\thDCc
\proclaim Lemma \thDCc. {\rm (cf.\ [\Jib], Proposition 2.20)} 
Let $\la,\mu\vdash n$. Then $\widetilde{G}^{(n)}_{\la\mu} 
\equiv\langle \la;\hw | \mu;\hw \rangle$ is given by 
\eqn\eqDCk{
\widetilde{G}^{(n)}_{\la\mu}  \eql \de_{\la\mu}\ b_\la(q)\,,
} 
where $b_\la(q) = \prod_{i\geq1} (q)_{m_i(\la)}$, and is independent
of $h$.\medskip

\proof Consider the reverse lexicographic ordering on the
set of partitions of $n$.
It is clear that $\langle \la;\hw | \mu;\hw\rangle = 0$ for
$\la > \mu$. Indeed, start by 
moving $T_{\la_1}$ to the right. This will produce terms 
with $T_m, m>\la_1$ and eventually give a vanishing contribution
when acting on the vacuum;  unless $T_{\la_1}$ 
is killed by a corresponding 
$T_{-\la_1}$, in which case we have  $\la_1=\mu_1$. And so on.
Because of symmetry, we  also have 
$\langle \la;\hw | \mu;\hw \rangle = 0$ for $\la < \mu$.
The diagonal terms are then easily calculated. The $h$-independence is
a consequence of the second simplification as discussed above. \Box\medskip

It immediately follows from Lemma \thDCc\ 
that the Kac determinant at level $n$ is given by
\eqn\eqDCl{
\widetilde{G}^{(n)} 
  \eql \prod_{\la\vdash n } \prod_{i\geq1} (q)_{m_i(\la)} 
  \eql \prod_{{\scr r,s\geq 1}\atop {\scr rs\leq n} }
  \left( 1-q^r \right)^{p(n-rs)}\,,
}
where we have used Lemma \thCBab.  One can verify that \eqDCl, indeed, 
agrees with the leading order behaviour of the Kac
determinant of Theorem \thCBa.
\thm\thDCe
\proclaim Theorem \thDCe. Let $q=\sqrtN{1}$. The vectors
\eqn\eqDCm{
|(n^N);\hw\rangle \eql (\widetilde{T}_{-n})^N |\hw\rangle\,,\qquad n\in\NN\,,
}
are singular for any $\hw\in\CC$.
\medskip

\proof From Lemma \thDCc\ it follows that,
for $q=\sqrtN{1}$, the vectors $(\wT_{-n})^N |\hw\rangle$ are
orthogonal to any vector in $\wM(\hw)$, i.e., they are in the radical
of $\wM(\hw)$, and hence  are primitive.%
\foot{See, e.g., [\BMP] for
the definition and properties of primitive vectors.}
However, from the explicit commutation relations \eqDCaa\ it easily follows
that $\wT_m (\wT_{-n})^N |\hw\rangle= 0$ for $m>0$.  Indeed,
commuting $\wT_m$ to the right will not produce
terms containing $\wT_{-p}$ for $p<n$, hence $\wT_m
(\wT_{-n})^N |\hw\rangle$ cannot be in the submodule
generated by the $(\wT_{-p})^N|\hw\rangle$ with $p<n$.  In
other words, it has to vanish. \Box\medskip

As expected, the existence of $h$-independent singular vectors is a
consequence of a large center of $\tvirq$, as described by the
following main theorem of this section.
\thm\thDCd
\proclaim Theorem \thDCd. For $q=\sqrtN{1}$, $N>2$, the center of 
$\tvirq$ is generated by the elements $(\widetilde{T}_n)^N, 
n\in\ZZ\backslash \{0\}$.\medskip

One can prove this theorem directly using the relations
\eqDCaa\ -- \eqDCad. We have summarized this rather tedious 
argument in Appendix~C. A more elegant proof 
will be given in due course using the free field
realization and symmetric functions.

The free field realization of $\tvirq$ is obtained by taking the
limit $t\to\infty$ in the free field realization of Section 2.3, after
rescaling the generators of ${\cH}_{p,q}$, 
\eqn\eqDresc{
q^{\be_0}\eql p^{-{1\over 2}}\,q^{\al_0}\,,\qquad
\be_m\eql-\al_m\, p^{{m\over 2}}\,,\qquad m\neq0\,.
}
In terms of the $q$-Heisenberg algebra $\widetilde{\cH}_q$,
\eqn\eqDCf{
[\be_m ,\be_n] \eql m\, (1-q^{|m|})\, \de_{m+n,0}\,,
}
we have a realization $\wimath: \tvirq \rightarrow 
\widetilde{\cH}_q$ defined by
\eqn\eqDCg{\eqalign{
\wimath(\wT_m) \eql \wLa^+_m\,,\qquad
&\wimath(\wT_{-m}) \eql \wLa^-_{-m}\,,\quad m>0\,,\cr
\wimath(\wT_0) &\eql \wLa^+_0 + \wLa^-_0\,,\cr}
}
where
\eqn\eqDCh{ \eqalign{
\wLa^+(z) & \eql \lim_{t\rightarrow\infty} \La^+( zp^{-{1\over 2}})
\eql q^{\be_0} \NO{\exp\left( \sum_{n\neq0} {\be_{n}\over n} z^{-n}
  \right)} \,,\cr
\wLa^-(z) &  \eql \lim_{t\rightarrow\infty} \La^-( zp^{{1\over 2}})
\eql q^{-\be_0} \NO{
 \exp\left( - \sum_{n\neq0} {\be_{n}\over n} z^{-n}
  \right)} \,.\cr} } 
The contraction identities now read (cf.\ \eqBBf)
\eqn\eqDCha{ 
\wLa^{\ep_1}(z_1) \wLa^{\ep_2}(z_2)  \eql 
\widetilde{f}^{\ep_1\ep_2}({z_2\over z_1}) 
\NO{\wLa^{\ep_1}(z_1) \wLa^{\ep_2}(z_2)}\,,
}
with
\eqn\eqBBg{
\widetilde{f}^{++}(x) \eql \widetilde{f}^{--}(x) \eql 
\widetilde{f}(x)^{-1}\,,\qquad \widetilde{f}^{+-}(x)\eql
\widetilde{f}^{-+}(x)\eql \widetilde{f}(x)\,,
}
where (cf.\ \eqDCda),
\eqn\eqDCi{
\widetilde{f}(x) \eql {1-qx\over 1-x}\,.
}

The modes of $\wLa^+(z)$ and  $\wLa^-(z)$ satisfy commutation
relations 
\eqnn\eqDCia\eqnn\eqDCxa\eqnn\eqDCxb
$$ \eqalignno{
[\wLa^+_m, \wLa^+_n]_q & \eql - [\wLa^+_{n+1},\wLa^+_{m-1}]_q \,,
&\eqDCia\cr
[\wLa^+_m, \wLa^-_n]_q & \eql - [\wLa^-_{n-1},\wLa^+_{m+1}]_q 
 + (1-q)^2\, \de_{m+n,0} \,,&\eqDCxa \cr
[\wLa^-_m, \wLa^-_n]_q & \eql - [\wLa^-_{n+1},\wLa^-_{m-1}]_q \,.
&\eqDCxb\cr}
$$
Note that \eqDCia\ and \eqDCxb\ are both equivalent to 
\eqDCab\ but with $m,n\in\ZZ$. 

\medskip
\ni {\bf Remark.} 
The algebra of the vertex operators $\wLa^\pm(z)$ has been studied in
[\Jib] (and in [\Jia] for $q=-1$) in the context of the
Hall-Littlewood symmetric functions.  In \eqDCf\ we have chosen a
slightly different normalization for $\wcH_q$, that is more suitable
to discuss the case $q=\sqrtN{1}$ than the one in [\Jib].
\medskip

Let $F(\be)$ be the Fock space of the Heisenberg algebra
$\widetilde{\cH}_q$. 
Since the structure of $\tvirq$ modules is
independent of $h$, we can set $\be= 0$ without loss of generality.
As in
Section 4.2 we identify $F(0)$ with the space of symmetric functions
$\La[q]$ by means of
\eqn\eqDCmc{
\widetilde\jmath\,(\be_{-\la_1}  \ldots \be_{-\la_\ell} |0\rangle) 
\eql p_\la(x)\,.
}
Note that $\widetilde\jmath$ is an isometry in the scalar
product on $\La[q]$ defined by \appAxca.

The specialization of Theorem \thDBc\ to the case $t\to\infty$ gives
an explicit identification of $\wimath(\widetilde M(h))$ with
symmetric functions (cf.\ [\Jib], Proposition 3.9).
\thm\thDCg
\proclaim Theorem \thDCg.  For any partition $\la$, we have 
\eqn\eqDCma{
\widetilde \jmath\,(\wLa^-_{-\la_1} \ldots \wLa^-_{-\la_n}
|0\rangle) \eql Q'_\la(x;q)\,, 
} 
where
\eqn\eqDCmb{
Q'_\la(x;q) \eql \Biggl( \prod_{i<j} {1-R_{ij} 
\over 1-qR_{ij}}\Biggr) \ h_{\la_1}(x)
\ldots h_{\la_n}(x)\,,
}
is a Milne polynomial.
\medskip

\proof  For a definition of Milne polynomials, see Appendix
A. The proof is essentially the same as in [\Jib] and can be found in
[\Gar]. \Box
\medskip

\ni
{\bf Remark.}
Note that \eqDCma\ makes sense for any sequence
$\la=(\la_1,\ldots,\la_n)$, where the $\la_i$ are nonnegative, but not
necessarily in descending order. Thus one can use \eqDCma\ to define
$Q'_\la(x;q)$ for such a general sequence $\la$. In Appendix A we have
introduced Milne polynomials for arbitrary sequences using their
relation to Hall-Littlewood polynomials. It is easy to see that the
two extensions are exactly the same, since the reordering identity
\appAxrs\ is identical with the commutation relations
\eqDCia\ (cf. [\MD], Example 2, p.\ 213). 
Moreover, \eqDCmb\ holds in this more general situation.
\medskip

A consequence of Theorem \thDCg\ is the following result, which
explicitly relates two types of ordering of 
products of modes of the vertex operators $\wLa^\pm(z)$.
\thm\thDCh
\proclaim Corollary \thDCh.  
\eqn\eqDCna{
\NO{ \wLa^-(z_1) \ldots \wLa^-(z_n) } ~\eql \sum_{\la_1\geq \la_2 \geq
\ldots \geq \la_n} \ P_\la(z;q)\, 
\wLa^-_{-\la_1} \ldots \wLa^-_{-\la_n} \,,
} 
where the sum is over all ordered sequences $\la$ (not necessarily positive).
\medskip

\proof Write 
\eqn\eqKj{
\NO{ \wLa^-(z_1) \ldots \wLa^-(z_n)} ~\eql \sum_{\la_1\geq \la_2 \geq
\ldots \geq \la_n} \ a_\la(z;q)\, 
\wLa^-_{-\la_1} \ldots \wLa^-_{-\la_n} \,.
}
Using \eqDCh\  we have 
\eqn\eqKk{ \eqalign{
\NO{ \wLa^-(z_1) \ldots \wLa^-(z_n) } |0\rangle
 &  \eql \exp\Biggl(
\sum_{k\geq1} \, {\be_{-k}\over k} (z_1^k + 
\ldots + z_n^k) \Biggr)|0\rangle \cr
&~\mathrel{\mathop=^{ \widetilde \jmath }}~
\exp\Biggl( \sum_{k\geq1}\, {1\over k}\, p_k(x) p_k(z)
 \Biggr) \eql \sum_\la \ P_\la(z;q) Q'_\la(x;q) \,,\cr}
}
where, in addition,  we have used
\eqn\eqKl{
\exp\Biggl( \sum_{k\geq1}\, {1\over k}\, p_k(x) p_k(z) \Biggr)
\eql \prod_{i,j} \left( {1\over 1-x_iz_j}\right) \,,
}
and the completeness relation \eqDCmh.

Now, on the one hand, Theorem \thDCg\ implies that, 
for partitions $\la_1\geq \la_2\geq \ldots\geq \la_n\geq0$
we have $a_\la(z;q) =P_\la(z;q)$.  On the other hand, from the `shift
invariance' of the commutators \eqDCia, it follows that
$a_{\la_1+k,\ldots,\la_n+k}(z;q) = 
(z_1\ldots z_n)^k a_{\la_1,\ldots,\la_n}(z;q)$ for all $k\in\ZZ$. This 
proves the corollary. \Box\medskip

We may now return to the main subject of this section, which is 
the center of $\tvirq$ and the proof of Theorem \thDCd.
\medskip

\ni {\it Proof of Theorem \thDCd:}  Let $q=\sqrtN{1}$. Consider the
expansion \eqDCna\ in Corollary \thDCh\ for $n=N$ and $z_i=z q^{i-1}$,
$i=1,\ldots,N$.  By the same argument as in the proof of Lemma~\thDBd,
we conclude that the coefficients $P_\la(z,zq,\ldots,zq^{N-1};q)$
vanish,  unless $m_i(\la)=N$ for some $i\in\ZZ_{\geq0}$, i.e.,
only the terms for which $\la=(m^N)$ for some $m\in\ZZ$ contribute to
the expansion.  Thus we find
\eqn\eqDCnc{
\NO{ \wLa^-(zq^{N-1})\ldots \wLa^-(zq)  \wLa^-(z) } ~\eql \sum_{m\in\ZZ} 
\ q^{{1\over2} mN(N-1)} (\wLa^-_{-m})^N \ z^{mN}\,.
}
Since $1+q^n+\ldots+q^{n(N-1)}=0$ for $n\not=0$
mod$\,N$, we also have
\eqn\exDCpfa{
\NO{ \wLa^-(zq^{N-1})\ldots \wLa^-(zq) \wLa^-(z) } ~\eql
\NO{ \exp\Biggl (\sum_{k\not=0}{\be_{-kN}\over k} z^{-kN}\Biggr)}\,.
} 
This shows that $(\wLa^-_{m})^N$, $m\in\ZZ$, are in the center of
$\widetilde{\cH}_q$. Since the free field realization acts faithfully
on $F(0)$, this also shows that (cf.\ \eqDCg) $(T_m)^N$, $m<0$, are
in the center of $\tvirq$, while repeating the same analysis for
$\wLa^+(z)$ extends this claim to $(T_m)^N$, $m>0$. 

One can check explicitly that, except for $N=2$, the $(T_0)^N$ are not
in the center, see the example in Appendix C. 

It remains to verify that we have identified the entire center. In
fact this follows from the irreducibility of the Verma module $M'(h)$,
which is proved by taking the $t\to\infty$ limit in  Theorem
\thDCfa, or directly from the Kac determinant in \eqDCl. 
This concludes the proof of the theorem. \Box\medskip

As a consequence of Theorems \thDCd\ and \thDCg, we obtain the
following interesting identity for Milne polynomials at $q=\sqrtN{1}$,
which was first discussed in [\LLTa,\LLTb] (see also [\MD], p.\ 234):
\thm\thDCi
\proclaim Corollary \thDCi.  Let $q=\sqrtN{1}$ and $n\in\NN$, then
\eqn\eqDCya{
Q'_{\la \cup (n^N)}(x;q) \eql Q'_{\la}(x;q) Q'_{(n^N)}(x;q)\,,
}
for any partition $\la$.
\medskip

\proof 
Suppose $\la_i\geq n$ for $i\leq k$ and $ \la_i<n$ for $i\geq
k+1$. Then, by  Theorem \thDCd,  we have
\eqn\eqCDdub{
\wLa^-_{-\la_1}\ldots \wLa^-_{-\la_k}(\wLa^-_{-n})^N
\wLa^-_{-\la_{k+1}}\ldots \wLa^-_{-\la_\ell}|0\rangle\eql
\wLa^-_{-\la_1}\ldots\wLa^-_{-\la_\ell} (\wLa^-_{-n})^N|0\rangle\,.
}
Since $(\wLa^-_{-n})^N$ is in the center of $\widetilde{\cH}_q$, the
action of $\wLa^-_{-\la_1}\ldots\wLa^-_{-\la_\ell}$ on the right hand
side is only through the creation operators, $\be_{-m}$, $m>0$, and
as a result $\widetilde \jmath$ factorizes, which proves \eqDCya. \Box
\medskip

Furthermore, note that in the context of the free field realization,
we obtain another simple proof of Theorem \thDCe\ 
by using Theorem \thDCg,
the faithfulness of $\widetilde{\imath}$ and the identity
\eqn\eqDCyb{\eqalign{
Q'_{(n^N)}(x;q) & \eql (-1)^{n(N-1)}\, h_n(x^N) \eql (-1)^{n(N-1)}\, 
 \sum_{\la\vdash n}\ z_\la^{-1}\, p_\la(x^N) \cr & \eql 
 (-1)^{n(N-1)}\, 
 \sum_{\la\vdash n}\ z_\la^{-1}\, p_{N\la}(x)\,,\cr}
}
which holds for $q=\sqrtN{1}$ (cf.\ [\MD], p.\ 235).

In Section 4.1 we have found that the center of $\virpq$ at $q=-1$
could be calculated in terms of symmetrized products of the
generators as given in \eqDAd. We will now show that a generalization
of this result holds for $\tvirq^\pm$ at $\qN$. In fact, this is a simple
consequence of the following symmetrization theorem which holds for an
arbitrary $q$.
\thm\thCDsymm
\proclaim Theorem \thCDsymm. For a  partition
$\la=(\la_1,\ldots,\la_n)$,
\eqn\eqCDsymm{
\sum_{\si\in S_n/S_n^\la} \widetilde T_{-\si\la_1}\ldots
\widetilde T_{-\si\la_n}\eql
\sum_{\mu} M_{\la\mu}(q)\bin{n}{m(\mu)}\widetilde T_{-\la_1}\ldots 
\widetilde T_{-\la_n}
\,,}
where $\si$ runs over all inequivalent permutations of the sequence
$(\la_1,\ldots,\la_n)$, $\mu$ over all partitions such that
$|\mu|=|\la|$ and $\ell(\mu)=\ell(\la)$, and the matrix $M(q)$ is given
by  $M(q)=K(1)^{-1}K(q)$, where $K(q)$ is the Kostka-Foulkes matrix.
\medskip

\proof Clearly, it is sufficient to verify \eqCDsymm\ when
both sides act on the vacuum in a Verma module. Then, by mapping the
Verma module to the symmetric functions, we find that in fact the
present theorem is equivalent to (the symmetrization) Lemma \thappAa\
for Milne polynomials in Appendix A. \Box\medskip

\ni
{\bf Remark.} There is an obvious counterpart of this result for
$\tvirq^+$. However, unlike in Section 4.1, there is no simple
symmetrization identity that would mix the positive and negative mode
generators of $\tvirq$. This might be one reason why there seems to be
no simple generalization of this theorem to  generic $t$.

\newsec{Discussion}

In Section 4.2 (see equation \eqDCnz)
we have seen that, for $q=\sqrtN{1}$, the character of 
the reduced Verma modules $M'(\hw)$ is given by
\eqn\eqGa{
\prod_{n\geq1}\left(\frac{1-x^{nN}}
{1-x^{n}}\right) \eql \prod_{n\geq1} (1 + x^n + \ldots + x^{n(N-1)})\,.
}
While the left hand side of this equation has a natural interpretation 
in terms of the character of a bosonic Fock space with the oscillators 
$\al_{nN},\, n\in\ZZ$, removed, the right hand side has a natural 
interpretation in terms of the character of the
Fock space of a so-called Gentile parafermion
of order $N-1$ [\Gen], 
i.e., a generalization of a fermion ($N=2$) defined by the
property that at most $N-1$ particles can occupy the same state, as
embodied in the equation $(\wT_n)^N=0,\, n\in\ZZ$, for $t\to\infty$ or
a deformation thereof for generic $t$.  In other words, 
the left hand side of this equation counts the number of partitions 
without parts of length $0\,{\rm mod}\,N$, while the right hand side counts 
the number of partitions with at most $N-1$ equal parts.  This equality 
is of course well-known in combinatorics.
 
Free fermions, of course, are well-known to give rise to conformal
field theories.  Indeed, as
we have noted in Section 4.2, we have a (non-canonical) action of
$\virpq$ at $q=\root 2\of{1} = -1$ on certain modules over the 
(undeformed) Virasoro algebra $\vir$ at central charge $c={1\over2}$.
This rather remarkable feature, that we have an action of a quantum 
group on modules over an undeformed algebra, 
resembles the action of Yangians on affine Lie algebra modules
(see [\BS] and references therein), and was one of our motivations 
to study the algebra $\virpq$ in the first place.
Of course, there is a crucial difference between between the two 
results, in that the affine Lie algebra modules are fully reducible 
into finite dimensional irreducible representations of this Yangian, while
the $c={1\over2}$ $\vir$ module carries an infinite dimensional 
irreducible representation of $\virpq$ at $q=-1$.  In fact, in Theorem
\thCc\ we have seen that $\virpq$ does not have any finite dimensional
irreducible representations in the category $\cO$ except for some trivial
ones.  It is an interesting open question whether $\virpq$ possesses 
any other nontrivial finite dimensional
irreducible representations (e.g., cyclic representations).

In addition, one may wonder whether, 
analogous to $q=-1$, the algebra $\virpq$ at $q=\sqrtN{1}$ can also
be realized on modules of $\vir$.
Clearly, the character \eqGa\ corresponds to the character of a conformal
field theory of central charge $c = (N-1)/N$ and is thus non-unitary
(for most $N$), as is a well-known property of Gentile parafermions 
for $N\neq2$ (see, e.g., [\Pol] and references therein).  This is left 
for further study as well.
Interestingly, a collection of $N$ 
such Gentile parafermions does realize a unitary conformal field theory,
namely $\widehat{\sln}$ at level $1$ [\Sch].

\medskip
\ni {\bf Acknowledgements:} P.B.\ is supported by a \qeii\ research 
fellowship from the Australian Research Council and K.P.\ is supported
in part by the U.S.\ Department of Energy Contract \#DE-FG03-84ER-40168.
P.B.\ would like to thank the University of Southern California for 
hospitality during the initial stages of this work and K.~Schoutens for 
discussions.\bigskip


\appendix{A}{Symmetric functions}

In this appendix we give basic definitions and summarize some results
from the theory of symmetric functions that are used throughout the
paper. For further details and omitted proofs the reader should
consult the monograph [\MD] and/or the references cited below.

\appsubsec{Symmetric functions}

Let $\La$ be the ring of symmetric functions in countably many
independent variables $x=\{x_1, x_2,\ldots\,\}$. Throughout the paper we are
using the following standard bases in $\La$ indexed by partitions,
$\la$:
\medskip

\item{(i)} monomial symmetric functions, $m_\la$,
\item{(ii)} elementary symmetric functions, $e_\la$,
\item{(iii)} complete symmetric functions, $h_\la$,
\item{(iv)} Schur functions, $s_\la$.
\medskip

One defines those functions as the limit, when the number of variables
goes to infinity, of the corresponding symmetric polynomials that are
stable with respect to the adjunction of variables. For instance, the
monomial symmetric function, $m_\la$, is the limit,
$m_\la(x)=\lim_{n\rightarrow\infty} m_\la(x_1,\ldots,x_n)$, of the
monomial symmetric polynomials
\eqn\appAxb{
m_\la(x_1,\ldots,x_n)=\sum x^{\al_1}_1\ldots x^{\al_n}_n\,,
}
where the sum in \appAxb\ runs over all inequivalent permutations
$\al=(\al_1,\ldots,\al_n)$ of the partition
$\la=(\la_1,\ldots,\la_n)$.

The elementary symmetric functions and the complete symmetric
functions are defined by
\eqn\appAxc{
e_\la=e_{\la_1}e_{\la_2}\ldots\,,\qquad h_\la=h_{\la_1}h_{\la_2}\ldots\,,
}
where $e_r$ and $h_r$ are determined from their generating series
\eqn\eqDBza{ \eqalign{
E(x;t) & \eql \sum_{n\geq0}\ e_n(x)\,t^n \eql \prod_{i\geq1}\, 
(1+x_it)\,,\cr
H(x;t) & \eql \sum_{n\geq0}\ h_n(x)\,t^n \eql \prod_{i\geq1}\, 
(1-x_it)^{-1}\,.\cr}
}

Finally, the Schur functions, $s_\la$, can be defined as  polynomials
in the elementary symmetric functions, $e_r$, or, equivalently, as
polynomials in the complete symmetric functions, $h_r$, namely
\eqn\appAxd{
s_\la=\det(e_{\la'_i-i+j})_{1\leq i,j\leq m}
\,,\qquad s_\la=\det(h_{\la_i-i+j})_{1\leq i,j\leq n}\,,
}
where $m\geq\ell(\la')$ and $n\geq \ell(\la)$,  respectively, and
$\la'$ is the partition conjugate to $\la$.

Upon extension of $\La$ to $\La_\QQ=\La\otimes_\ZZ\QQ$, the ring of
symmetric functions with rational coefficients, one can introduce yet
another basis, namely the power sum symmetric functions,
$p_\la$. Those are defined by
\eqn\appAxe{
p_\la \eql p_{\la_1}p_{\la_2}\ldots \,,\qquad
p_r\eql\sum x_i^r=m_{(r)}\,.
}

The standard scalar product on $\La_\QQ$ is defined by
\eqn\appAxf{
\langle p_\la,p_\mu\rangle \eql z_\la\de_{\la\mu}\,,
} where $z_\la$ is given in \eqCBab. In fact \appAxf\ induces a
well-defined scalar product on $\La$, with respect to which the Schur
functions, $s_\la$, form an orthonormal basis.

\appsubsec{Hall-Littlewood and Milne symmetric functions}

The Hall-Littlewood (HL) symmetric functions, $P_\la(x;q)$, form an
orthogonal basis in the space of one parameter symmetric functions
$\La[q]=\La\otimes_\ZZ \ZZ[q]$ with the scalar product defined by
\eqn\appAxca{
\langle p_\la,p_\mu\rangle_q \eql z_\la\,\de_{\la\mu}\,
\prod_{i=1}^{\ell(\la)} (1-q^{\la_i})\,.
}  
They can be calculated by applying the
Gramm-Schmidt orthogonalization algorithm to the basis of Schur
functions and are given explicitly as the
limit of the corresponding HL polynomials [\Lit]
\eqn\appAxff{
P_\la(x_1,\ldots,x_n;q)\eql {1\over v_\la(q)} 
\sum_{w\in S_n} w\Biggl(
x_1^{\la_1}\ldots x_n^{\la_n} \prod_{i<j}{x_i-qx_j\over x_i-x_j}\Biggr)\,, 
}
where
\eqn\appAxg{
v_m(q)\eql{(q)_m\over (1-q)^m}\,, \qquad v_\la(q)\eql\prod_{i\geq
0}v_{m_i}(q)\,.
}  
It is understood in \appAxff\ that the permutations $w$ act on the
variables $x_1,\ldots,x_n$.

The functions $P_\la(x;q)$ interpolate between the Schur
functions, $s_\la(x)$, and the monomial symmetric functions,
$m_\la(x)$, namely
\eqn\appAxh{
P_\la(x;0)\eql s_\la(x)\,,\qquad P_\la(x;1)\eql m_\la(x)\,.
}

The transition matrix $K(q)$ defined by
\eqn\appAxxv{
s_\la(x)\eql\sum_\mu K_{\la\mu}(q)\,P_\mu(x;q)\,, } is strictly upper
unitriangular with respect to the natural order of partitions, i.e.,
$K_{\la\mu}(q)=0$ unless $|\la|=|\mu|$ and $\la\succeq\mu$, and
$K_{\la\la}(q)=1$. The entries $K_{\la\mu}(q)\in\ZZ[q]$ are called the
Kostka-Foulkes polynomials. One can show  that their coefficients
are non-negative integers. It follows from \appAxh\ that $K(1)$ is the
transition matrix between the monomial symmetric functions and the
Schur functions,
\eqn\appAzx{
s_\la(x)=\sum_{\mu}K_{\la\mu}(1)\,m_\mu(x)\,.}
Its entries are called the Kostka numbers.

Another family of symmetric functions, $Q_\la(x;q)$, also referred to as HL
symmetric functions, are scalar multiples of the $P_\la(x;q)$, defined
by
\eqn\appAxz{
Q_\la(x;q)\eql b_\la(q)\,P_\la(x,q)\,,  } 
where $b_\la(q) = \prod_{i\geq1} (q)_{m_i(\la)}$, 
so that
\eqn\appAxpr{
\langle P_\la,Q_\mu\rangle_q=\de_{\la\mu}\,.}
An obvious consequence of the definition \appAxz\ (see, also \appAyc)
is that for $\qN$ the functions $Q_{(k^N)}(x;q)$ vanish.

It has already been observed in [\Lit] that the definition of the
$Q_\la(x;q)$ in \appAxff\ and \appAxz\ makes sense for any sequence of
integers $(\la_1,\ldots,\la_n)$ that are not necessarily in a
descending order and/or are not positive. For such generalized $Q_\la$
one can prove [\Lit]  a reordering identity that allows to reduce
$Q_\la$, when $\la_i$ are not in descending order, to a linear
combination of the $Q_\mu$, where $\mu_i$ are in descending order. For
a two-term sequence the formula is (see, [\MD] p.\ 214)
\eqn\appAxrs{
Q_{(m,n)}-q\, Q_{(n,m)}\eql  
-Q_{(n-1,m+1)}+ q\, Q_{(m+1,n-1)}\,.
}
The same equality holds within sequences $\la$ of length greater than two.

Now, let us consider the symmetric polynomials corresponding to
$Q_\la$. It is clear that a particularly simple polynomial arises if
we set the number of variables equal to the length of the partition,
namely
\eqn\appAya{
Q_\la(x_1,\ldots,x_n;q)\eql (1-q)^n \sum_{w\in S_n} w\Biggl(
x_1^{\la_1}\ldots x_n^{\la_n} \prod_{i<j}{x_i-qx_j\over
x_i-x_j}\Biggr)\,, } where $\ell(\la)=n$.  The main result of this
section is a symmetrization lemma for such polynomials (cf. [\Mob],
Theorem 1):

\thm\thappAa
\proclaim Lemma \thappAa. Let
$\la=(\la_1,\ldots,\la_n)$ be a partition and $S_n/S_n^\la$ the 
subgroup of distinct permutations of the sequence $\la$. Then
\eqn\KPscoj{
\sum_{\si\in S_n/S_n^\la}Q_{\si\la}(x_1,\ldots,x_n;q)
=\sum_{\mu}M_{\la\mu}(q)\bin{n}{m(\mu)}Q_\mu(x_1,\ldots,x_n;q)\,,
}
where $M(q)=K(1)^{-1}K(q)$, and the sum on the r.h.s.\ in \KPscoj\
runs only over permutations $\mu$ such that $|\mu|=|\la|$ and
$\ell(\mu)=\ell(\la)$.
\medskip

\proof By applying the symmetrization in $\la_1,\ldots,\la_n$ to the
right hand side of \appAya, we obtain
\eqn\appAyc{
\eqalign{
\sum_{\si\in S_n/S_n^\la}Q_{\si\la}(x_1,\ldots,x_n;q)&\eql
(1-q)^n\sum_{\si\in S_n/S_n^\la}
\sum_{w\in S_n}w\Biggl(x_1^{\si\la_1}\ldots 
x_n^{\si\la_n}\prod_{i<j}{x_i-qx_j\over x_i-x_j}\Biggr)\cr
&\eql
(1-q)^n\sum_{\si\in S_n/S_n^\la}x_1^{\si\la_1}\ldots 
x_n^{\si\la_n}
\sum_{w\in S_n}w\Biggl(\prod_{i<j}{x_i-qx_j\over x_i-x_j}\Biggr)\cr
&\eql
(q)_nm_\la(x_1,\ldots,x_n)\,,\cr}
}
where we used \appAxb\ and the identity (see, [\MD] p.\ 207)
\eqn\appAxin{
\sum_{w\in S_n}w \Biggl(\prod_{i<j}{x_i-qx_j\over
x_i-x_j}\Biggr)=v_n(q)\,.
}
Now use the fact that 
\eqn\appApfl{
m_\la(x_1,\ldots,x_n) = m_\la(x_1,\ldots,x_m)\big|_{ x_{n+1}=\ldots=x_m=0}\,,
}
and
\eqn\appAflp{
Q_\mu(x_1,\ldots,x_m)\big|_{ x_{n+1}=\ldots=x_m=0} = \cases{
  Q_\mu(x_1,\ldots,x_n) & if $\ell(\mu) \leq n\,,$ \cr
  0 & if $\ell(\mu) > n\,,$\cr}
}
together with \appAxxv, \appAzx\ and \appAxz\
to find
\eqn\eqqqqone{
m_\la(x_1,\ldots,x_n) = \sum_{\mu}\  
M_{\la\mu}(q)\ {1\over b_\mu(q)}\ 
Q_\mu(x_1,\ldots,x_n;q)\,,
}
where 
\eqn\appAeqqq{
M(q) = K(1)^{-1} K(q)\,.  } The restriction on the range of the sum
over $\mu$ is then a straightforward consequence of the strict upper
unitriangularity of $K(q)$. \Box
\medskip

We also need to introduce another family of symmetric functions,
$Q'_\la(x;q)$, called Milne's symmetric functions [\Mila,\Milb] (see,
also [\Gar,\LLTa,\LLTb]).  They are related to the HL functions,
$Q_\la(x;q)$, by a change of variables 
\eqn\eqDCmf{
Q'_\la(x;q) \eql Q_\la({x\over 1-q};q)\,,
}
which has to be understood the sense of the $\la$-ring notation. This
means that $Q'_\la(x;q)$ is the image of $Q_\la(x;q)$ by the
ring homomorphism of $\La[q]$ that sends $p_r(x)$ to $(1-q^r)^{-1}
p_r(x)$. 
A defining property of Milne's functions is their orthogonality to
the HL functions,
\eqn\appAmo{
\langle P_\la,Q'_\mu\rangle\eql \de_{\la\mu}\,,}
with respect to the product \appAxf. This is equivalent to the
completeness relation 
\eqn\eqDCmh{
\sum_\la \ P_\la(x;q)\, Q'_\la(y;q) \eql \prod_{i,j}\left( {1\over 1-x_iy_j}
\right)\,.
}
Finally, an expansion in terms of Schur functions yields
\eqn\eqDCmg{
Q'_\la(x;q) \eql \sum_{\mu } \ K_{\mu\la}(q)\, s_\mu(x)\,.
}

One can use \eqDCmf\ to extend the definition of $Q'_\la(x;q)$ to
arbitrary sequences $\la$. Given that the transformation defined in
\eqDCmf\ is a ring homomorphism, it is clear that so defined
$Q'_\la(x;q)$ will satisfy both the reordering identity \appAxrs\ as
well as \KPscoj\ of the symmetrization lemma.


\appendix{B}{Explicit singular vectors for $q=\sqrtN{1}$}

In Theorem \thDBb\ of Section 4.2 we have seen that
\eqn\eqAPAaa{
\Psi(z) \eql \lim_{z_i\rightarrow zq^{N-i}}\ 
\Biggl( \prod_{i<j}f({z_j\over z_i})\Biggr) \,
T(z_1)\ldots T(z_N) |\hw\rangle\,,
}
is a (well-defined) generating series of singular vectors in $M(\hw)$
for $q=\sqrtN{1}$.  We have also outlined how to make sense out of this
expression.  
Here we carry out the procedure for $N=3$.  First,
\thm\thAPAb
\proclaim Theorem \thAPAb.  We have 
\eqn\eqAPAca{ \eqalign{
\Biggl( \prod_{i<j} & f( {z_j\over z_i} ) \Biggr)
T(z_1) T(z_2) T(z_3)|\hw\rangle \cr & \eql 
\sum_{\scr m_1,m_2,m_3\geq 0} z_1^{m_1}z_2^{m_2}z_3^{m_3} 
\Biggl( \prod_{i<j} f(R_{ij})\Biggr) \, T_{-m_1}T_{-m_2}T_{-m_3}
|\hw\rangle \cr
& \qquad + \ze
\Big({ p {(z_3/ z_2)} \over 1-p {(z_3/ z_2)} }
F^{-+}({z_3\over z_1})-{ p^{-1} {(z_3/ z_2)} \over 1-p^{-1} {(z_3/ z_2)}}
F^{+-}({z_3\over z_1})\Big) T(z_1)|\hw\rangle\cr
&\qquad + \ze \Big(
{p{(z_3/ z_1)}\over 1-p{(z_3/ z_1)}}
F_{(d)}^{-+}({z_3\over z_2})-{p^{-1}{(z_3/ z_1)}\over 1-p^{-1}{(z_3/ z_1)}}
F_{(d)}^{+-} ({z_3\over z_2})\Big)T(z_2)|\hw\rangle\cr
&\qquad +\ze  \Big({1\over 1-p{(z_2/ z_1)}}F^{+-}(p{z_3\over z_1})-
{1\over 1-p^{-1}{(z_2/ z_1)}}
F^{-+}(p^{-1}{z_3\over z_1})\Big)T(z_3)|\hw\rangle\,.\cr}
}
where
\eqn\eqAPAcb{
F_{(d)}^{\ep_1\ep_2}(x) \eql \sum_{0\leq m\leq d}\ F^{\ep_1\ep_2}_m \, x^m\,.
}
\medskip

We sketch the proof, which is based on the following
\thm\thAPAa
\proclaim Lemma \thAPAa.  Let 
\eqn\eqAPAza{
c_{m,n} \eql \ze\, (p^m - p^{-m})\, \de_{m+n,0}\,,
}
then we have
\eqn\eqAPAzb{ \eqalign{
\sum_{l_1,l_2\geq0}\ f_{l_1} f_{l_2} T_{-m_1-l_1-l_2} 
c_{-m_2+l_1,-m_3+l_2} & \eql
\ze \, \left( p^{-m_2} F^{-+}_{m_2+m_3}-
p^{m_2}F^{+-}_{m_2+m_3} \right) T_{-m_1-m_2-m_3}\,,\cr
\sum_{l_1,l_2\geq0}\ f_{l_1} f_{l_2} c_{-m_1-l_1,-m_2-l_2} 
T_{-m_3+l_1+l_2} & \eql
\ze \,\left( p^{-m_1}F^{+-}_{-m_1-m_2}-p^{m_1}F^{-+}_{-m_1-m_2}
\right)T_{-m_1-m_2-m_3}\,,\cr}
}
where $F^{\pm\mp}(x)=\sum_m F^{\pm\mp}_m x^m$ is defined in \eqBBi.\medskip

\proof We have
\eqn\eqKm{ \eqalign{
\sum_{l_1,l_2\geq0}\ & f_{l_1} f_{l_2} 
T_{-m_1-l_1-l_2} c_{-m_2+l_1,-m_3+l_2}  \cr & \eql
\ze \sum_{p,q\geq0}\ f_{l_1} f_{l_2} 
( p^{-m_2+l_1} - p^{m_2-l_1}) \de_{-m_2-m_3+l_1+l_2}
T_{-m_1-l_1-l_2} \cr
& \eql \ze \left ( p^{-m_2} f(px)f(x){\big|_{x^{m_2+m_3}} }- p^{m_2} 
f(p^{-1}x)f(x){\big|_{x^{m_2+m_3}}}\right) T_{-m_1-m_2-m_3}\,.\cr}
}
The other identity is shown similarly. \Box\medskip

\ni {\it Proof of Theorem \thAPAb:} Writing out the left hand side of 
\eqAPAca\ in modes we distinguish three cases\hfil\break
\line{\hfill
(i) $m_1\geq0, m_2\geq0, m_3 \geq 0$,\quad
(ii) $m_2<0$, $m_3\geq 0$, \quad (iii) $m_1<0$, $m_2\geq 0$, $m_3\geq 0$.
\hfill}
In case (i) only a finite number of terms contribute.  This gives the first
term on the right hand side of \eqAPAca.  In case (ii) we move 
the term $T_{-m_2+l_{12}-l_{23}}$ to the right using the commutator \eqBAa.
Then use Lemma \thAPAa\ to find 
$$
\ze\, \sum_{m_1}\sum_{m_2\leq -1}\sum_{m_3\geq 0}\
z_1^{m_1}z_2^{m_2}z_3^{m_3} 
\left (p^{-m_2}F^{-+}_{m_2+m_3}-
p^{m_2}F^{+-}_{m_2+m_3} \right) T_{-m_1-m_2-m_3} |\hw\rangle\,.
$$
The sum over $m_1$ is unrestricted. The nonvanishing terms must have
$m_2+m_3\geq 0$ because of the moding of $F^{+-}$ and $F^{-+}$. Thus,
with the restriction on $m_2$ in place, we may let the sum on $m_3$ 
run over all integers.  This allows for a change in the summation
\eqn\eqKn{
m'=m_2+m_3\,,\qquad m''=m_1+m_2+m_3\,,
}
after wich all sums can be performed and we obtain 
$$
\ze
\Big({ p {(z_3/ z_2)} \over 1-p {(z_3/ z_2)} }
F^{-+}({z_3\over z_1})-{ p^{-1} {(z_3/ z_2)} \over 1-p^{-1} {(z_3/ z_2)}}
F^{+-}({z_3\over z_1})\Big) T(z_1)|\hw\rangle\,.
$$
The remaining case (iii) is analyzed similarly. \Box\medskip

By taking the limit $z_i\to zq^{N-i}$ in Theorem \thAPAb\ and 
using $F^{-+}(q^{-1}) = F^{-+}(p^{-1}q^{-2}) = F^{+-}(q^{-2}) = 0$,
we obtain
\thm\thAPAc
\proclaim Corollary \thAPAc.  The following is an explicit form for the 
level $d$ singular vector at $q=\root3\of1$ and arbitrary $\hw\in\CC$
\eqn\eqAPAba{ \eqalign{
\Psi_d \eql & \sum_{\scr m_1,m_2,m_3\geq 0\atop\scr
m_1+m_2+m_3=d} q^{2m_1+m_2} 
\Biggl( \prod_{i<j} f(R_{ij})\Biggr) T_{-m_1}T_{-m_2}T_{-m_3}|\hw\rangle\cr
&+ \ze 
\Big[ q^{2d}\Big({1\over qp^{-1} -1}F^{-+}(q^{-2})\Big)
T_{-d}\cr
&+ q^d\Big({1\over q^2p^{-1} -1}F_{(d)}^{-+}(q^{-1})-{1\over q^2p-1}
F_{(d)}^{+-}(q^{-1})\Big)T_{-d}\cr
&+ \Big({qp^{-1}\over qp^{-1}-1}F^{+-}(pq^{-2})\Big)T_{-d}\Big]
|\hw\rangle\,.\cr}
}
\medskip

The expression $\Psi_d$ vanishes for $d\neq 0\,{\rm mod}\,3$ (we have
explicitly verified this for small values of $d$).  

{}From the structure of $f_{l}$, and the duality
symmetry $(p,q)\to (p^{-1},q^{-1})$, it follows that the singular 
vectors at $q=\sqrtN{1}$ for 
$d=0\,{\rm mod}\,N$ can be written in terms of the 
following set of invariants
\eqn\eqAPAaa{
\De^{rs}_{m_1\ldots m_k}=a_{rs}{p^rq^s(1+p^{n-2r}q^{N-2s})
\over \prod_{i=1}^k(1+p^i)^{m_i}}\,,\qquad n=\sum_i im_i\,,
}
where the normalization factor 
\eqn\eqAPAab{
a_{rs}=\cases{ {1\over 2} & for $2r=n$ and $2s=0\ \mod N\,,$\cr
		1 & otherwise$\,.$\cr}
}
Using \eqAPAba\ we now find the 
following singular vectors at $q=\root 3 \of 1$ for arbitrary
$\hw\in\CC$.  At level $d=3$
\eqn\eqAPAa{
\Psi_3 \eql \left( T_{-1}T_{-1}T_{-1} + a_{210}\, T_{-2}T_{-1}T_{0} + a_{300}\,
T_{-3}T_{0}T_{0} + a_3\, T_{-3} \right) |\hw\rangle\,,
}
with
\eqn\eqAPAb{ \eqalign{
a_{210} & \eql  -3\De^{10}_{2}\,,\cr
a_{300} & \eql  -3\De^{21}_{31}\,,\cr
a_{3} & \eql 3 \De_{11}^{11} \,,\cr}
}
and, at level $d=6$
\eqn\eqAPAc{ \eqalign{
\Psi_6 \eql & \Big(T_{-2}T_{-2}T_{-2} + a_{321}\,T_{-3}T_{-2}T_{-1} 
+ a_{411}\,T_{-4}T_{-1}T_{-1} +
a_{330}\, T_{-3}T_{-3}T_{0}  \cr & + a_{420}\, T_{-4}T_{-2}T_{0} + 
a_{510}\, T_{-5}T_{-1} T_{0} + a_{600}\, T_{-6}T_{0}T_{0} + a_{6}\, T_{-6} 
\Big) 
|\hw\rangle\,,\cr}
}
where
\eqn\eqAPAd{ \eqalign{
a_{321} & \eql -3\De^{10}_{2}\,,\cr
a_{411} & \eql -3\De^{21}_{31}\,,\cr
a_{330} & \eql -3\De^{21}_{31}\,,\cr
a_{420} & \eql -3\De^{40}_{42}+3\De^{30}_{41}-3\De^{21}_{31}\,,\cr
a_{510} & \eql -6\De^{40}_{42}+3\De^{31}_{52}\,,\cr
a_{600} & \eql 27\De^{80}_{6301}+12\De^{60}_{63}+18\De^{70}_{6201}+
6\De^{50}_{62}-3\De^{61}_{5301}-3\De^{51}_{5301}-3\De^{41}_{5301}\,,\cr
a_{6}   & \eql -3\De^{60}_{4201}+3\De^{20}_{4201}-6\De^{40}_{42}
-12\De^{50}_{4101}+6\De^{30}_{41}-3\De^{51}_{4201}+3\De^{41}_{4201}
+6\De^{31}_{4201}\,. \cr}
}
The equality of the coefficients $a_{210}=a_{321}$ and $a_{300}=a_{411}$
in the singular vectors at $d=3$ and $d=6$ is a direct consequence of
Lemma \thDBd\ (ii) and the expression \eqDBj.

In addition we have computed the level $d=4$ singular vector at 
$q=\root 4 \of 1$
\eqn\eqAPAe{ \eqalign{
\Psi_4 \eql & \Big( T_{-1}T_{-1}T_{-1}T_{-1} + 
a_{2110}\, T_{-2}T_{-1}T_{-1}T_{0} + 
a_{2200}\, T_{-2}T_{-2}T_{0}T_{0} +
a_{3100}\, T_{-3}T_{-1}T_{0}T_{0} \cr & +
a_{4000}\, T_{-4}T_{0}T_{0}T_{0} +
a_{22}\, T_{-2}T_{-2}+ a_{31}\, 
T_{-3}T_{-1} + a_{40}\, T_{-4}T_{0}\Big) |\hw\rangle\,,\cr}
}
where
\eqn\eqAPAf{ \eqalign{
a_{2110} & \eql -4 \De^{10}_2 - 4 \De^{11}_3 \,,\cr
a_{2200} & \eql 8 \De^{20}_4 \,,\cr
a_{3100} & \eql 12 \De^{20}_4 + 4 \De^{21}_5 \,,\cr
a_{4000} & \eql -8 \De^{30}_{501} -8 \De^{31}_{501}\,,\cr
a_{22}   & \eql -8 \De^{10}_{2}   \,,\cr
a_{31}   & \eql -4 \De^{10}_{2} + 4 \De^{11}_{3}\,,\cr
a_{40}   & \eql 8\De^{20}_{301} - 8 \De^{21}_{301} \,.\cr}
}

Observe that for $t\to\infty$ we have $\De^{r,s}_{m_1\ldots m_k}
\to 0$ for $0<r\leq 2n$.  Thus, in all the examples above, in the
limit $t\to\infty$ only the leading term $(T_{-n})^N|\hw\rangle$
survives.  Furthermore, note that, 
both for $N=3$ and $N=4$, all the factors in front of the
fundamental invariants $\De^{rs}_{m_1\ldots m_k}$ are a multiple of $N$
(except for the one of the leading term).  Thus, it appears that
calculating modulo $N$, in the appropriate sense, is somehow equivalent 
to considering the $t\to\infty$ limit.

\appendix{C}{The center revisited}

\def\virinf{\tvirq}
\def\TT{t}
\def\SS{s}

In this appendix we establish in a direct way some elementary
identities for the products of generators of $\tvirq^\pm$ for a
generic $q$.  Specialization of those results to $\qN$ yields a direct
proof of Theorem \thDCa.

\appsubsec{Preliminaries}
The defining relations \eqDCab\ of $\virinf^\pm$ are invariant under
the rescaling $\widetilde T_m\rightarrow a^{|m|} \widetilde T_m$,
$a\in\CC$, while the remaining relations in \eqDCaa\ -- \eqDCad\ can be
further simplified by a judicious choice of $a$. In the following we
will find it convenient to work with the generators
\eqn\eqKo{
\TT_m \eql q^{{|m|\over 2}}\, \wT_m \,,\qquad m\in\ZZ\,.
}
Then the relations involving the generators of $\virinf^+$ are given
by
\eqn\KPcmtwo{
\TT_{ m}\TT_{ n} \eql q 
\TT_{ n}\TT_{ m}
-(1-q ) \sum_{l=1}^{m-n-1} \TT_{ m-l}
\TT_{ n+l}\,,\qquad m>n\geq 1\,,
}
\eqn\KPcmfour{
\TT_{ m}\TT_0 \eql q \TT_0\TT_{ m}
-(1-q^{})\sum_{l=1}^{m-1}\TT_{ m-l}
\TT_l+(q-q^{-1})\sum_{l=1}^\infty\TT_{-l}\TT_{m+l}\,,
\qquad m\geq 1\,,
}
\eqn\KPcmone{
\TT_{m}\TT_{-n} \eql q\TT_{-n}\TT_m+
(q-q^{-1})\sum_{l=1}^\infty
\TT_{-n-l}\TT_{m+l}+(1-q)q^{m}\de_{m,n}\,,\qquad m,n\geq 1\,.
}
while the remaining ones have a similar form and can easily be worked
out.

\appsubsec{Identities in $\virinf^\pm$}

Now let us consider the subalgebras $\virinf^\pm$ in more detail. We
will discuss explicitly only the case of $\virinf^+$ with the
extension to $\virinf^-$ being obvious.

The relation \KPcmtwo\ implies that two subsequent generators
satisfy
\eqn\appXa{
\TT_{m+1}\TT_m \eql q\TT_m\TT_{m+1}\,,
}
For the generators $\TT_{m+k}$ and $\TT_m$ with $k\geq 2$, there are
additional  terms, though it is
still possible to rewrite \KPcmtwo\ in a symmetric form
\eqn\appXb{
\sum_{j=1}^k \TT_{m+j}\TT_{m+k-j} \eql
q\sum_{j=1}^k\TT_{m+k-j}\TT_{m+j}\,.
}

In view of \appXb\ it is then natural to consider the
following sums of generators%
\foot{Since \KPcmtwo\ is invariant under rescaling, there seems to be
no advantage in working with the generating series $t^+(z)=\sum_{n\geq
1}t_{ n} z^{-n}$. }
\eqn\eqKp{
\SS_m \eql \TT_m+\TT_{m+1}+\ldots\,,
}
in term of which \KPcmtwo\ is equivalent to
\eqn\KPCls{
\SS_{m+1}\TT_m \eql q\TT_m\SS_{m+1}-(1-q)\SS_{m+1}\SS_{m+1}\,.
}
The infinite sums and their products here and below should be
understood in the graded sense.  By iterating \KPCls\ we prove

\thm\KPcrlem
\proclaim Lemma \KPcrlem. For $n\geq 1$,
\eqn\KPCgenla{
(\SS_{m+1})^n\TT_m \eql q^n\TT_m(\SS_{m+1})^n-(1-q^n)(\SS_{m+1})^{n+1}\,,
}
\eqn\KPCgenlb{
\SS_{m+1}(\TT_m)^n~=~\sum_{j=0}^n (-1)^jq^{n-j}(q)_j
	\bin{n}{j}(\TT_m)^{n-j}(\SS_{m+1})^{j+1}\,.
}

\medskip
\noindent
Then,  a straightforward induction yields

\thm\KPcmlma
\proclaim Lemma \KPcmlma. For $n\geq 1$,
\eqn\KPCprlem{
(\TT_m+\SS_{m+1})^n \eql
\sum_{j=0}^n q^{{1\over 2}j(j-1)} \bin{n}{j}(\TT_m)^{n-j}
(\SS_{m+1})^j\,,
}
\eqn\KPcmlsaa{
\left(q\TT_m+(q-q^{-1})\SS_{m+1}\right)^n \eql \sum_{j=0}^n
(-1)^j{q^{n-2j}(q)_{j+1}\over 1-q}\bin{n}{j}(\TT_m)^{n-j}(\SS_{m+1})^j\,.
}
\medskip

For a partition $\la$ of length $\ell(\la)=n$, let
$\la^{{\rm op}}=(\la^{{\rm op}}_1,\ldots,\la^{{\rm op}}_n)$ 
be the increasing sequence of positive
integers, $\la^{{\rm op}}_i=\la_{n-i}$, $i=1,\ldots,n$. 
For such a sequence we define
\eqn\KPhtdef{\eqalign{
\height(\la^{{\rm op}})&\eql\sum_{i=1}^{\ell(\la)} (\ell(\la)-i)
(\la^{{\rm op}}_i-1)\cr
&\eql n(\la)-{1\over 2}\ell(\la)(\ell(\la)-1)\,.\cr}
}
Finally, in terms of multiplicities, we 
can write $(\la_1,\ldots,\la_n)=(1^{m_1}2^{m_2}\ldots\,)$, where
$m_i=m_i(\la)$.

Using Lemma~\KPcmlma\ we will now establish the main result of this
section, which gives the expansion of powers of $\SS_m$ into ordered
products of generators.
\thm\KPsymth
\proclaim Theorem \KPsymth. For $n\geq 1$,
\eqn\KPCcnsum{
(\SS_m)^n \eql
\sum_{\{\la\,|\, \ell(\la)=n\}}q^{{\height} (\la^{{\rm op}})}
\bin{n}{m(\la)}
\TT_{m+\la^{{\rm op}}_1-1}\ldots \TT_{m+\la^{{\rm op}}_n-1}\,.
}
\medskip

\proof 
We write $\SS_m=\TT_m+\SS_{m+1}$, and expand
$(\SS_m)^n=(\TT_m+\SS_{m+1})^n$ using \KPcmlsaa,
\eqn\eqKq{
(\SS_m)^n \eql \sum_{m_1+n_2=n}q^{\height(1^{m_1}2^{n_2})}
\bin{n}{m_1,n_2}(\TT_{m})^{m_1}(\SS_{m+1})^{n_2}\,.
}
Since at a given level only products of a finite number of generators,
$\TT_m,\TT_{m+1},\ldots,\TT_{m+s-1}$, can appear, after repeating this
expansion $s-1$ times we obtain a multiple sum 
\eqn\eqKr{
(\SS_m)^n \eql \sum_{m_1+\ldots+m_s=n}\!\!\!
q^{\sum_{i=1}^{s-1}\height(1^{m_i}2^{n_{i+1}})}
\bin{n}{m_1,n_2}\ldots
\bin{n_{s-1}}{m_{s-1},m_s}
(\TT_m)^{m_1}\ldots (\TT_{m+s-1})^{m_s}\,,
}
where $n_i=m_i+\ldots+m_s$. It follows from definition \KPhtdef\
that
\eqn\eqKs{
\sum_{i=1}^{s-1}{\height}(1^{m_i}2^{{n_{i+1}}}) \eql {\height}(1^{m_1}
2^{m_2}\ldots\, s^{m_s})\,,
}
and it is obvious that
\eqn\eqKt{
\bin{n}{m_1,n_2}\bin{n_2}{m_2,n_3}\ldots 
\bin{n_{k-1}}{m_{s-1},m_s}\eql \bin{n}{m_1,\ldots,m_s}\,.
}
This completes the proof of the theorem. \Box

\appsubsec{Proof of Theorem \thDCd}

\thm\KPcent
\proclaim Theorem \KPcent. Let $\qN$. Then $(\TT_m)^N$, $m\geq 1$, 
belong to the center of\/ $\virinf$.
\medskip

\proof We will show separately that 
\eqn\KPtoshow{
(\TT_m)^N\TT_n \eql \TT_n(\TT_m)^N\,,\qquad m\geq 1\,,
}
for $\pm n>0$ and $n=0$.
\smallskip

\noindent
Case 1: $n>0$

{}From Lemma \KPcrlem\ we find
\eqn\eqKu{
(\SS_{m+1})^N\TT_m \eql \TT_m(\SS_{m+1})^N\,,
}
and
\eqn\eqKv{
\SS_{m+1}(\TT_m)^N \eql (\TT_m)^N\SS_{m+1}\,,
}
for arbitrary $m\geq 1$. Since, by Theorem~\KPsymth
\eqn\eqKw{
(\SS_m)^N \eql \sum_{k=0}^\infty q^{{1\over 2}kN(N-1)}(\TT_{m+k})^N\,, 
}
where each term in the sum is at a different level, this 
implies \KPtoshow. \smallskip

\noindent
Case 2: $-n>0$

First rewrite \KPcmone\ in a more convenient form,
\eqn\KPcmab{
\TT_m\TT_{-n} \eql \sum_{l=1}^\infty
a_l \TT_{-n-l}\TT_{m+l}+c_n\de_{m,n}\,,
}
where $a_0=q$, $a_l=q-q^{-1}$, $l\geq 1$, and $c_n=(1-q)q^{-n}$.  By 
repeated use of \KPcmab\ we obtain
\eqn\eqKx{
\eqalign{
(\TT_m)^N\TT_{-n}&\eql c_m(\TT_m)^{N-1}\de_{m,n}\cr&\quad+
\sum_{k=1}^{N-1}(\TT_m)^{N-k-1}\sum_{l_1,\ldots,l_k=0}^\infty
a_{l_1}\ldots a_{l_k}c_{n+l_1+\ldots+l_k}\TT_{m+l_1}\ldots \TT_{m+l_k}
\de_{m,n+l_1+\ldots+l_k}\cr&\quad+
\sum_{l_1,\ldots,l_N=0}^\infty 
a_{l_1}\ldots a_{l_N}\TT_{-n-l_1-\ldots-l_N}
\TT_{m+l_1}\ldots \TT_{m+l_N}\,.\cr}
}
This may be recognized as
\eqn\KPcmgenb{\eqalign{
(\TT_m)^N\TT_{-n}
&\eql c_m(\TT_m)^{N-1}\de_{m,n}\cr &\quad +c_m\sum_{k=1}^{N-1}
(\TT_m)^{N-k-1}\big(\, q\TT_m+(q-q^{-1})\SS_{m+1}\,\big)^k_{(k+1)m-n}\cr
&\quad+
\sum_{l\geq 0}\TT_{-n-l}\big(\, q\TT_m+(q-q^{-1})\SS_{m+1}\,\big)^N_{Nm+l}\,,
\cr}
}
where the subscripts on the brackets indicate the level.

For $m<n$, we find using \KPcmlsaa\ that only the last term in
\KPcmgenb\ with $l=0$ contributes giving $\TT_{-n}(\TT_m)^N$. 

For $m=n$, all terms in \KPcmgenb\ contribute, however, we may set
$\SS_{m+1}=0$ in the second term. Thus the central charge term has an
overall factor of $1+q+\ldots+q^{N-1}=0$, while the last term is the
same as above.

For $m>n$, the first term does not contribute. Let us
first consider the second term proportional to the central charge. Using
\KPcmlsaa\ we can rewrite it as
$$
c_m
\sum_{k=1}^{N-1}\sum_{j=0}^k (-1)^j{q^{k-2j}(q)_{j+1}\over 1-q}
\bin{k}{j}(\TT_m)^{N-j-1}(\SS_{m+1})^j\big|_{Nm-n}\,,
$$
where only the terms at the level $Nm-n$ contribute. In particular, at
this level we must have $j\geq 1$. Changing the order of summation we obtain 
\eqn\eqKy{ \eqalign{
&c_m\sum_{j=1}^{N-1}
(-1)^j{q^{-j}(q)_{j+1}\over 1-q}\big(\, \sum_{k=j}^{N-1}q^{k-j}\bin{k}{j}
\,\big) (\TT_m)^{N-j-1}(\SS_{m+1})^j
\big|_{Nm-n}\cr
&\eql c_m
\sum_{j=1}^{N-1}
(-1)^j{q^{-j}(q)_{j+1}\over 1-q}\bin{N}{j+1}
(\TT_m)^{N-j-1}(\SS_{m+1})^j
\big|_{Nm-n}\,.\cr}
}
Clearly all terms in the sum vanish. Thus once more the only
contribution arises from the last term in \KPcmgenb\ and it is the same
as above.
\smallskip

\noindent
Case 3: $n=0$

Using \KPcmfour\ we find
\eqn\eqKz{
\eqalign{
(\TT_m)^n\TT_0 &\eql q^N\TT_0(\TT_m)^N\cr
&\quad + \sum_{k=1}^N q^{k-1}(\TT_m)^{N-k}[\,(q-1)\sum_{l=1}^{m-1}
\TT_{m-l}\TT_l\,](\TT_m)^{k-1}\cr
&\quad + \sum_{k=1}^N q^{k-1}(\TT_m)^{N-k}[\,(q-q^{-1})\sum_{l=1}^\infty
\TT_{-l}\TT_{m+l}\,](\TT_m)^{k-1}\,.\cr}
}
This expression is clearly the same as the one obtained by setting
$\TT_0=\TT_0^++\TT_0^-$, where $\TT_0^+$ satisfies \KPcmtwo\ with
$n=0$ and $\TT_0^-$ satisfies \KPcmone\ with $n=0$. Since the proof of
the two cases above does not depend on the specific value of $n$, by
exactly the same algebra we show that $\TT_0$ commutes with
$(\TT_m)^N$. \Box\medskip

An obvious modification of the above argument proves that $(T_m)^N$
for $m<0$ lies in the center. This  then concludes the proof of
Theorem~\thDCd. \medskip

\ni {\bf Remark.}   Note that $(\TT_0)^N$, ($N\geq3$), is not in the
center as the following example for $N=3$ shows,
\eqn\eqKaa{ \eqalign{
(\TT_0)^3\TT_{-2}&\eql q^3\TT_{-2}(\TT_0)^3-(1-q)q^2(1+q+q^2)
(\TT_{-1})^2(\TT_0)^2\cr
&\quad +(1-q)^3q(2+q)\TT_{-2}\TT_0-(1-q)^4(1+q)
(\TT_{-1})^2 +\ldots\,,\cr}
}
where the dots stand for terms with strictly positive modes on the right.

\footatend\immediate\closeout\rfile
\baselineskip=14pt{\bigskip\noindent {\bf References}}%
\bigskip{\frenchspacing%
\parindent=20pt\escapechar=` \input refs.tmp\vfill\eject}\nonfrenchspacing
\vfil\eject\end